\input amstex
\documentstyle{amsppt}
\magnification=\magstep 1
\document
\vsize6.7in

\chardef\oldatsign=\catcode`\@
\catcode`\@=11
\newif\ifdraftmode			
\global\draftmodefalse



%
\font@\twelverm=cmr12 
\font@\twelvei=cmmi12 \skewchar\twelvei='177 
\font@\twelvesy=cmsy10 scaled\magstep1 \skewchar\twelvesy='060 
\font@\twelveex=cmex10 scaled\magstep1 
\font@\twelvemsa=msam10 scaled\magstep1 
\font@\twelvemsb=msbm10 scaled\magstep1 
\font@\twelvebf=cmbx12 
\font@\twelvett=cmtt12 
\font@\twelvesl=cmsl12 
\font@\twelveit=cmti12 
\font@\twelvesmc=cmcsc10 scaled\magstep1 
%
%
\font@\ninerm=cmr9 
\font@\ninei=cmmi9 \skewchar\ninei='177 
\font@\ninesy=cmsy9 \skewchar\ninesy='60 
\font@\ninemsa=msam9
\font@\ninemsb=msbm9
\font@\ninebf=cmbx9
%
%
%
\font@\ttlrm=cmbx12 scaled \magstep2 
\font@\ttlsy=cmsy10 scaled \magstep3 
\font@\tensmc=cmcsc10 
%
%
\def\normaltype{
	\def\pointsize@{12}%
	\abovedisplayskip18\p@ plus5\p@ minus9\p@
	\belowdisplayskip18\p@ plus5\p@ minus9\p@
	\abovedisplayshortskip1\p@ plus3\p@
	\belowdisplayshortskip9\p@ plus3\p@ minus4\p@
	\textonlyfont@\rm\twelverm
	\textonlyfont@\it\twelveit
	\textonlyfont@\sl\twelvesl
	\textonlyfont@\bf\twelvebf
	\textonlyfont@\smc\twelvesmc
	\ifsyntax@
		\def\big##1{{\hbox{$\left##1\right.$}}}%
	\else
		\let\big\twelvebig@
 \textfont0=\twelverm \scriptfont0=\ninerm \scriptscriptfont0=\sevenrm
 \textfont1=\twelvei  \scriptfont1=\ninei  \scriptscriptfont1=\seveni
 \textfont2=\twelvesy \scriptfont2=\ninesy \scriptscriptfont2=\sevensy
 \textfont3=\twelveex \scriptfont3=\twelveex  \scriptscriptfont3=\twelveex
 \textfont\itfam=\twelveit \def\it{\fam\itfam\twelveit}%
 \textfont\slfam=\twelvesl \def\sl{\fam\slfam\twelvesl}%
 \textfont\bffam=\twelvebf \def\bf{\fam\bffam\twelvebf}%
 \scriptfont\bffam=\ninebf \scriptscriptfont\bffam=\sevenbf
 \textfont\ttfam=\twelvett \def\tt{\fam\ttfam\twelvett}%
 \textfont\msafam=\twelvemsa \scriptfont\msafam=\ninemsa
 \scriptscriptfont\msafam=\sevenmsa
 \textfont\msbfam=\twelvemsb \scriptfont\msbfam=\ninemsb
 \scriptscriptfont\msbfam=\sevenmsb
	\fi
 \normalbaselineskip=\twelvebaselineskip
 \setbox\strutbox=\hbox{\vrule height12\p@ depth6\p@
      width0\p@}%
 \normalbaselines\rm \ex@=.2326ex%
}
%
%
%
\def\smalltype{
	\def\pointsize@{10}%
	\abovedisplayskip12\p@ plus3\p@ minus9\p@
	\belowdisplayskip12\p@ plus3\p@ minus9\p@
	\abovedisplayshortskip\z@ plus3\p@
	\belowdisplayshortskip7\p@ plus3\p@ minus4\p@
	\textonlyfont@\rm\tenrm
	\textonlyfont@\it\tenit
	\textonlyfont@\sl\tensl
	\textonlyfont@\bf\tenbf
	\textonlyfont@\smc\tensmc
	\ifsyntax@
		\def\big##1{{\hbox{$\left##1\right.$}}}%
	\else
		\let\big\tenbig@
	\textfont0=\tenrm \scriptfont0=\sevenrm \scriptscriptfont0=\fiverm 
	\textfont1=\teni  \scriptfont1=\seveni  \scriptscriptfont1=\fivei
	\textfont2=\tensy \scriptfont2=\sevensy \scriptscriptfont2=\fivesy 
	\textfont3=\tenex \scriptfont3=\tenex \scriptscriptfont3=\tenex
	\textfont\itfam=\tenit \def\it{\fam\itfam\tenit}%
	\textfont\slfam=\tensl \def\sl{\fam\slfam\tensl}%
	\textfont\bffam=\tenbf \def\bf{\fam\bffam\tenbf}%
	\scriptfont\bffam=\sevenbf \scriptscriptfont\bffam=\fivebf
	\textfont\msafam=\tenmsa
	\scriptfont\msafam=\sevenmsa
	\scriptscriptfont\msafam=\fivemsa
	\textfont\msbfam=\tenmsb
	\scriptfont\msbfam=\sevenmsb
	\scriptscriptfont\msbfam=\fivemsb
		\textfont\ttfam=\tentt \def\tt{\fam\ttfam\tentt}%
	\fi
 \normalbaselineskip 14\p@
 \setbox\strutbox=\hbox{\vrule height10\p@ depth4\p@ width0\p@}%
 \normalbaselines\rm \ex@=.2326ex%
}

\def\titletype{
	\def\pointsize@{17}%
	\textonlyfont@\rm\ttlrm
	\ifsyntax@
		\def\big##1{{\hbox{$\left##1\right.$}}}%
	\else
		\let\big\twelvebig@
		\textfont0=\ttlrm \scriptfont0=\twelverm
		\scriptscriptfont0=\tenrm
		\textfont2=\ttlsy \scriptfont2=\twelvesy
		\scriptscriptfont2=\tensy
	\fi
	\normalbaselineskip 25\p@
	\setbox\strutbox=\hbox{\vrule height17\p@ depth8\p@ width0\p@}%
	\normalbaselines
	\rm
	\ex@=.2326ex%
}

\def\tenbig@#1{
	{%
		\hbox{%
			$%
			\left
			#1%
			\vbox to8.5\p@{}%
			\right.%
			\n@space
			$%
		}%
	}%
}

\def\twelvebig@#1{%
	{%
		\hbox{%
			$%
			\left
			#1%
			\vbox to10.2\p@{}
			\right.%
			\n@space
			$%
		}%
	}%
}

%
%
%
%
%
\newif\ifl@beloutopen
\newwrite\l@belout
\newread\l@belin

\global\let\currentfile=\jobname

\def\getfile#1{%
	\immediate\closeout\l@belout
	\global\l@beloutopenfalse
	\gdef\currentfile{#1}%
	\input #1%
	\par
	\newpage
}

\def\getxrefs#1{%
	\bgroup
		\def\gobble##1{}
		\edef\list@{#1,}%
		\def\gr@boff##1,##2\end{
			\openin\l@belin=##1.xref
			\ifeof\l@belin
			\else
				\closein\l@belin
				\input ##1.xref
			\fi
			\def\list@{##2}%
			\ifx\list@\empty
				\let\next=\gobble
			\else
				\let\next=\gr@boff
			\fi
			\expandafter\next\list@\end
		}%
		\expandafter\gr@boff\list@\end
	\egroup
}

\def\testdefined#1#2#3{%
	\expandafter\ifx
	\csname #1\endcsname
	\relax
	#3%
	\else #2\fi
}

\def\document{%
	\minaw@11.11128\ex@ 
	\def\alloclist@{\empty}%
	\def\fontlist@{\empty}%
	\openin\l@belin=\jobname.xref	
	\ifeof\l@belin\else
		\closein\l@belin
		\input \jobname.xref
	\fi
}

\def\getst@te#1#2{%
	\edef\st@te{\csname #1s!#2\endcsname}%
	\expandafter\ifx\st@te\relax
		\def\st@te{0}%
	\fi
}

\def\setst@te#1#2#3{%
	\expandafter
	\gdef\csname #1s!#2\endcsname{#3}%
}

\outer\def\setupautolabel#1#2{%
	\def\newcount@{\global\alloc@0\count\countdef\insc@unt}	
	\def\newtoks@{\global\alloc@5\toks\toksdef\@cclvi}
	\expandafter\newcount@\csname #1Number\endcsname
	\expandafter\global\csname #1Number\endcsname=1%
	\expandafter\newtoks@\csname #1l@bel\endcsname
	\expandafter\global\csname #1l@bel\endcsname={#2}%
}

\def\reflabel#1#2{%
	\testdefined{#1l@bel}
	{
		\getst@te{#1}{#2}%
		\ifcase\st@te
			???
			\message{Unresolved forward reference to
				label #2. Use another pass.}%
		\or	
			\setst@te{#1}{#2}2
			\csname #1l!#2\endcsname 
		\or	
			\csname #1l!#2\endcsname 
		\or	
			\csname #1l!#2\endcsname 
		\fi
	}{
		{\escapechar=-1 
		\errmessage{You haven't done a
			\string\\setupautolabel\space for type #1!}%
		}%
	}%
}

{\catcode`\{=12 \catcode`\}=12
	\catcode`\[=1 \catcode`\]=2
	\xdef\Lbrace[{]
	\xdef\Rbrace[}]%
]%

\def\setlabel#1#2{%
	\testdefined{#1l@bel}
	{
		\edef\templ@bel@{\expandafter\the
			\csname #1l@bel\endcsname}%
		\def\@rgtwo{#2}%
		\ifx\@rgtwo\empty
		\else
			\ifl@beloutopen\else
				\immediate\openout\l@belout=\currentfile.xref
				\global\l@beloutopentrue
			\fi
			\getst@te{#1}{#2}%
			\ifcase\st@te
			\or	
			\or	
				\edef\oldnumber@{\csname #1l!#2\endcsname}%
				\edef\newnumber@{\templ@bel@}%
				\ifx\newnumber@\oldnumber@
				\else
					\message{A forward reference to label 
						#2 has been resolved
						incorrectly.  Use another
						pass.}%
				\fi
			\or	
				\errmessage{Same label #2 used in two
					\string\setlabel s!}%
			\fi
			\expandafter\xdef\csname #1l!#2\endcsname
				{\templ@bel@}
			\setst@te{#1}{#2}3%
			\immediate\write\l@belout 
				{\string\expandafter\string\gdef
				\string\csname\space #1l!#2%
				\string\endcsname
				\Lbrace\templ@bel@\Rbrace
				}%
			\immediate\write\l@belout 
				{\string\expandafter\string\gdef
				\string\csname\space #1s!#2%
				\string\endcsname
				\Lbrace 1\Rbrace
				}%
		\fi
		\templ@bel@	
		\expandafter\ifx\envir@end\endref 
			\gdef\marginalhook@{\marginal{#2}}%
		\else
			\marginal{#2}
		\fi
		\expandafter\global\expandafter\advance	
			\csname #1Number\endcsname
			by 1 %
	}{
		{\escapechar=-1
		\errmessage{You haven't done a \string\\setupautolabel\space
			for type #1!}%
		}%
	}%
}


\newcount\SectionNumber
\setupautolabel{t}{\number\SectionNumber.\number\tNumber}
\setupautolabel{r}{\number\rNumber}
\setupautolabel{T}{\number\TNumber}

\define\rref{\reflabel{r}}
\define\tref{\reflabel{t}}

\define\tnum{\setlabel{t}}
\define\rnum{\setlabel{r}}

%
\def\strutdepth{\dp\strutbox}%
\def\strutheight{\ht\strutbox}%

\newif\iftagmode
\tagmodefalse

\let\old@tagform@=\tagform@
\def\tagform@{\tagmodetrue\old@tagform@}

\def\marginal#1{%
	\ifvmode
	\else
		\strut
	\fi
	\ifdraftmode
		\ifmmode
			\ifinner
				\let\Vorvadjust=\Vadjust
			\else
				\let\Vorvadjust=\vadjust
			\fi
		\else
			\let\Vorvadjust=\Vadjust
		\fi
		\iftagmode	
			\llap{%
				\smalltype
				\vtop to 0pt{%
					\pretolerance=2000
					\tolerance=5000
					\raggedright
					\hsize=.72in
					\parindent=0pt
					\strut
					#1%
					\vss
				}%
				\kern.08in
				\iftagsleft@
				\else
					\kern\hsize
				\fi
			}%
		\else
			\Vorvadjust{%
				\kern-\strutdepth 
				{%
					\smalltype
					\kern-\strutheight 
					\llap{%
						\vtop to 0pt{%
							\kern0pt
							\pretolerance=2000
							\tolerance=5000
							\raggedright
							\hsize=.5in
							\parindent=0pt
							\strut
							#1%
							\vss
						}%
						\kern.08in
					}%
					\kern\strutheight
				}%
				\kern\strutdepth
			}
		\fi
	\fi
}


\newbox\Vadjustbox

\def\Vadjust#1{
	\global\setbox\Vadjustbox=\vbox{#1}%
	\ifmmode
		\ifinner
			\innerVadjust
		\fi		
	\else
		\innerVadjust
	\fi
}

\def\innerVadjust{%
	\def\nexti{\aftergroup\innerVadjust}%
	\def\nextii{%
		\ifvmode
			\hrule height 0pt 
			\box\Vadjustbox
		\else
			\vadjust{\box\Vadjustbox}%
		\fi
	}%
	\ifinner
		\let\next=\nexti
	\else
		\let\next=\nextii
	\fi
	\next
}%

\global\let\marginalhook@\empty

\def\endref{%
\setbox\tw@\box\thr@@
\makerefbox?\thr@@{\endgraf\egroup}%
  \endref@
  \endgraf
  \endgroup
  \keyhook@
  \marginalhook@
  \global\let\keyhook@\empty 
  \global\let\marginalhook@\empty 
}

\catcode`\@=\oldatsign

\nologo

\def\depth{\operatorname{depth}}
\define\iso{\cong}
\define \Hgy{\operatorname{H}}
\def\incl{\operatorname{incl}}
\redefine\b{\beta}
\define\link{\:\!}
\def\a{\alpha}
\def \f {\varphi} 
\def \p {\oplus}
\def \t {\otimes}

\def\and{\quad\text{and}\quad}
\def\blank{\underline{\phantom{X}}}
\define\({\left(}
\define\){\right)}
\define\[{\left[}
\define\]{\right]}
\define \g{\gamma}
\define \onto {\twoheadrightarrow}

\define\cok{\operatorname{coker}}
\define\W{\tsize\bigwedge}
\def \Sym{\operatorname{Sym}}
\define\grade{\operatorname{grade}}
\define\w{\wedge}
\define\htt{\operatorname{ht}}
\define \maxm {\frak m}
\define\tW{\W}
\def\reg{\operatorname{reg}}
\def\F{\Bbb F}
\def\GL{\operatorname{GL}}
\def\m{\maxm}
\topmatter
\abstract
 Consider a height two ideal,
  $I$,    which is minimally generated by $m$ homogeneous forms of degree $d$ in the polynomial ring  $R=k[x,y]$. Suppose   that one column in the homogeneous presenting matrix $\f$ of $I$ has  entries of degree $n$ and all of the other entries of $\f$ are linear. We identify an explicit generating set for the ideal $\Cal A$ which defines the Rees algebra $\Cal R=R[It]$; so $\Cal R=S/\Cal A$ for the polynomial ring $S=R[T_1,\dots,T_m]$. We resolve $\Cal R$ as an $S$-module and $I^s$ as an $R$-module, for all powers $s$. The proof uses the homogeneous coordinate ring, $A=S/H$, of  a  rational normal scroll,   with $H\subseteq \Cal A$. The ideal $\Cal AA$ is isomorphic to the $n^{\text{th}}$ symbolic power of a height one prime ideal $K$ of $A$. The ideal $K^{(n)}$ is generated by monomials. Whenever possible, we study $A/K^{(n)}$ in place of $A/\Cal AA$ because the generators of $K^{(n)}$ are much less complicated then  the generators of $\Cal AA$. We obtain a filtration of $K^{(n)}$ in which the factors are polynomial rings, hypersurface rings, or modules resolved by generalized Eagon-Northcott complexes. The generators of $I$ parameterize an algebraic curve $\Cal C$ in projective $m-1$ space. The defining equations of the special fiber ring $\Cal R/(x,y)\Cal R$ yield a solution of the implicitization problem for $\Cal C$. 
\endabstract

\title Rational Normal Scrolls and the Defining Equations of Rees Algebras \endtitle
  \leftheadtext{Kustin, Polini, and Ulrich}
\rightheadtext{Defining Equations of Rees Algebras}
 \author Andrew R. Kustin\footnote{Supported in part by the National Security Agency.\hphantom{the National Science Foundation and XXX}}, Claudia Polini\footnote{Supported in part by the National Science Foundation and the National Security Agency.\hphantom{XXX}}, and Bernd  Ulrich\footnote{Supported in part by the National Science Foundation.\hphantom{and the National Security Agency XXX}}\endauthor
 \address
Mathematics Department,
University of South Carolina,
Columbia, SC 29208\endaddress
\email kustin\@math.sc.edu \endemail
\address Mathematics Department,
University of Notre Dame,
Notre Dame, IN 46556\endaddress 
\email     cpolini\@nd.edu \endemail
\address 
Department of Mathematics, 
Purdue University,
West Lafayette, IN 47907
\endaddress
\email          ulrich\@math.purdue.edu \endemail
\endtopmatter

\document


\bigpagebreak

\SectionNumber=0\tNumber=1

\flushpar{\bf   Introduction.}\footnote""{\hskip-.17in 2000 {\it Mathematics Subject Classification.} Primary: 14Q05, Secondary: 13C20, 13D02, 14E05, 14M12, 68W30.}\footnote""{\hskip-.17in {\it Key words and phrases.} Divisorial ideal, Eagon-Northcott complex,  Hilbert-Burch Theorem, Implicitization,  Rational normal scroll,  Rees algebra, Special fiber ring.}

\medskip

In this paper we address the problem of determining the equations
that define the Rees algebra of an ideal. Recall that the {\it Rees
algebra \/} $\Cal{R}(I)$ of an ideal $I$ in a commutative ring $R$
is the graded subalgebra $R[It]$ of the polynomial ring $R[t]$. Any
finite generating sequence $\delta_1,  \ldots, \delta_m$ of $I$ gives rise to
an $R$-algebra epimorphism
$$\define\trightarrow{@>>\pretend\haswidth{\text{Cliff}}>} \Pi: R[T_1, \ldots , T_m]  \trightarrow \Cal{R}(I)
$$ mapping $T_i$ to $\delta_it$, whose kernel is the {\it defining ideal
\/} $\Cal{A}$ of the Rees algebra. Besides encoding asymptotic
properties of the powers of an ideal, the Rees algebra realizes,
algebraically, the blow-up of a variety along a subvariety. Though
blowing up is a fundamental operation in the birational study of
algebraic varieties and, in particular, in the process of
desingularization, an explicit description of the resulting variety
in terms of defining equations remains a difficult problem. In other
words, the structure and shape of the ideal $\Cal{A}$ defining the
Rees algebra is still wide open, though it has been the subject of a
long list of articles over the past thirty years that is too
extensive to quote here. Much of this work requires that the Rees
algebra have the Cohen-Macaulay property, which allows for better
control over the number and the degrees of the defining equations.
In this paper instead, we will discard of this assumption and deal
with a class of ideals whose Rees algebras are never Cohen-Macaulay.

We consider the case where $R=k[x,y]$ is a polynomial ring over a
field $k$ and $I$ is an ideal of height two minimally generated by
$m$ forms $\delta_1, \ldots, \delta_m$ of degree $d$. The Hilbert Burch
Theorem guarantees that $I$ is generated by the maximal order minors
of an $m$ by $m-1$ matrix $\f$ with homogeneous entries of constant
degree along each column. Thus, in addition to $m$ and $d$ the other
important piece of data are the column degrees of $\f$. In the
present paper, the column degrees of $\f$ are $(1,\dots,1,n)$. In
other words, the entries of one column of $\f$ have arbitrary degree
$n$, all of the other entries of $\f$ are linear; we say that the
ideal $I$ is {\it almost linearly presented}. In this setting we are
able to identify homogeneous generators of the defining ideal $\Cal
A$ of the Rees ring $\Cal{R}(I)$. We can safely assume that $n \geq
2$, for otherwise $I=(x,y)^d$ and the answer is well known (see, for
instance, \cite{\rref{MU}}). Incidentally, except when $n=1$, the Rees ring
$\Cal{R}(I)$ is never Cohen-Macaulay. Hong, Simis and Vasconcelos
\cite{\rref{HSV}} had identified the ideal $\Cal{A}$ if $m=3$ and $n \leq
5$, and they proposed a conjectural, inductive procedure for finding
a generating set of $\Cal{A}$ if $n$ is arbitrary. Their conjecture
was proved in \cite{\rref{CHW}}, thus solving the case of arbitrary almost
linearly presented almost complete intersection ideals in two
variables. Whereas the method of \cite{\rref{HSV}} and \cite{\rref{CHW}} is based
on iterations of `Jacobian duals' and `Sylvester determinants', our
approach is entirely different and allows for closed formulas for
all defining equations at once, besides avoiding the need to
restrict the number of generators of $I$.

To determine the defining ideal $\Cal{A}$ of the Rees ring one often
uses the fact that its presentation map $\Pi$ factors through the
symmetric algebra $\text{Sym}(I)$. It then remains to determine the
kernel of the natural epimorphism
$$\text{Sym}(I) \longrightarrow \Cal{R}(I) \, ,$$
since the defining ideal of $\text{Sym}(I)$ can be described easily.
On the downside however, $\text{Sym}(I)$ does not have good
ring-theoretic properties in general, for instance, it is hardly
ever a domain. Thus, the main idea in our approach is to replace
$\text{Sym}(I)$ by a different algebra $A$ that is still `closer' to
$\Cal{R}(I)$ than the polynomial ring $R[T_1, \ldots, T_m]$, but is
a normal domain unlike the symmetric algebra $\text{Sym}(I)$. We
prove that in our setting, the algebra $A$ can be chosen to be the
homogeneous coordinate ring of a three-dimensional rational normal
scroll. The map $\Pi$ induces an epimorphism
$$A \longrightarrow \Cal{R}(I) \, ,$$
whose kernel is a height one prime ideal $\Cal{A}A$ of the normal
domain $A$, and hence gives rise to an element of the divisor class
group group of $A$. Now in \cite{\rref{KPU-d}} we study divisors on
rational normal scrolls of arbitrary dimension -- most notably, for
any given divisor class we describe an explicit monomial generating
set of an unmixed ideal representing it, and we investigate these
ideals in detail by providing free resolutions. Exhibiting an
isomorphism between $\Cal{A}A$ and the much simpler, monomial
representative of its divisor class we obtain closed formulas for
the
defining equations of $\Cal{R}(I)$ (Theorem \tref{A=L}), which turn
out to be tremendously complicated despite the seemingly strong
assumptions on $I$! We go on to compute the depth of $\Cal{R}(I)$
(Theorem \tref{T1}), the reduction number of $I$ (Theorem
\tref{red1}), the Hilbert function and Castelnuovo-Mumford
regularity of all powers $I^s$ (Corollary \tref{C2.6}), and even the
minimal resolution of $I^s$ (Theorem \tref{pwrs4}). To do so, we
replace the ideal $\Cal{A}$ by its simpler, isomorphic model, which
is the $n$th symbolic power $K^{(n)}$ of a prime ideal $K$ generated
by the images of linear forms in $A$.  We then construct a
filtration of $K^{(n)}$ whose factors are easy to study; these
factors turn out to be Eagon-Northcott modules -- in the sense that
they are resolved by generalized Eagon-Northcott complexes. Our
answers are very different depending on whether the linear part of
the matrix $\f$ has a generalized row of zeros, meaning a row of
zeros after elementary row operations.

Finding the defining ideal $\Cal{A}$ of the Rees ring solves, in
particular, another classical problem in elimination theory: An
ideal $I \subset R=k[x,y]$ of height two generated by forms $\delta_1,
\ldots, \delta_m$ of degree $d$ gives rise to a morphism $\Bbb P^1
\longrightarrow \Bbb P^{m-1}$ mapping $[x_0:y_0]$ to
$[\delta_1(x_0,y_0):\cdots :\delta_m(x_0,y_0)]$. The image of this map is a
curve $\Cal C \subset \Bbb P^{m-1}$ with homogeneous coordinate ring
$k[\delta_1, \ldots, \delta_m]$. The latter ring is isomorphic to the {\it
special fiber ring \/}  $\Cal{F}(I)$ of $I$, which is defined as $k
\otimes_R \Cal{R}(I) $ with $x$ and $y$ acting trivially on $k$.
Thus we pay attention to the depth and algebraic properties of the
ring $\Cal{F}(I)$ as well (Theorem \tref{T1}). Clearly, the defining
ideal $\Cal{A}$ of the Rees ring yields, in particular, the defining
ideal $\Cal{A}k[T_1, \ldots, T_m]$ of the special fiber ring and
hence of the curve $\Cal{C}$. The problem of finding the defining
ideal of a curve $\Cal{C}$ that is given parametrically has also
attracted the attention of the geometric modeling community, where
it is known as `implicitization problem'; see, for example,
\cite{\rref{B01},\rref{BC05},\rref{BJ03},\rref{EM04}}. One technique
used in solving this problem is the method of `moving curves' and
the `moving curve ideal', which is nothing but the defining ideal
$\Cal{A}$ of the Rees ring \cite{\rref{C08},\rref{C06},\rref{B07}}.

This paper is organized as follows. In Section one we establish the
connection between the Rees algebra $\Cal R(I)$ and a rational
normal scroll ring. 
In Section two we find a canonical form for the linear part of the
matrix which presents the ideal $I$. The form is used in Section
one; furthermore, this form allows us to calculate the regularity of
$I^s$ for all $s$. We record an explicit generating set for the
defining ideal of  $\Cal R(I)$ in Section three. In Section four we
calculate the reduction number of $I$, the regularity and depth of
$\Cal F(I)$, and the depth  $\Cal R(I)$. The filtration of $K^{(n)}$
by Eagon-Northcott modules is in Section five. In Section six we
resolve $I^s$ and verify the regularity calculation of Section two.

\smallskip

The graded ring $S=\bigoplus\limits_{i\ge 0}S_i$ is a {\it standard graded $S_0$-algebra} if $S$ is generated as an $S_0$-algebra by $S_1$ and $S_1$ is finitely generated as an $S_0$-module. In this discussion $R$ is a standard graded polynomial ring over a field.  If 
$N$ is a finitely generated non-zero graded $R$-module and 
$$0\to F_k\to \dots \to F_0\to N\to 0,$$ with $F_{i}=\bigoplus_{j=1}^{\b_i}R(-t_{i,j})$, is the minimal homogeneous resolution of $N$ by free $R$-modules, then 
the  Castelnuovo-Mumford {\it regularity}  of $N$  is equal to  $\reg(N) =\max_{i,j} \{t_{i,j} - i\}$.

 \definition{Convention}Throughout this paper, $k$ is a  field; every ring $A$ that we consider is graded and finitely generated as an algebra over  $A_0=k$; and every $A$-module $M=\bigoplus M_i$ that we consider is graded and finitely generated. We use $\lambda_A(\blank)$ for the length of an $A$-module. It follows that 
$$\lambda_k(M_i),\quad \lambda_A\(\frac{\bigoplus_{i\le j}M_j}{\bigoplus_{i< j}M_j}\),\and \dim_k(M_i)$$ are equal. We write $\lambda(M_i)$ for the common value. Of course, all three numbers are equal to the value of the Hilbert function $\Hgy_M$ at $i$, denoted $\Hgy_M(i)$.
\enddefinition

  \definition{Convention}  For each  statement  ``S'', we define  $$\chi(\text{S})= \cases 1,&\text{if S is true, and}\\ 0,&\text{if S is false.}\endcases$$In particular, $\chi(i=j)$ has the same value as the Kronecker delta $\delta_{ij}$. \enddefinition

\definition{Notation} If $\theta$ is a real number, then $\lceil \theta\rceil$ and $\lfloor \theta \rfloor$ are the ``round up'' and ``round down'' of $\theta$, respectively; that is, $\lceil \theta\rceil$ and $\lfloor \theta \rfloor$ are the integers with
$$ \lceil \theta \rceil-1<\theta \le \lceil \theta \rceil\and
\lfloor \theta \rfloor\le \theta < \lfloor \theta \rfloor+1.
$$\enddefinition

For any terms or concepts that we neglected to define, consult \cite{\rref{BH}} first.
\bigskip
\bigskip
\bigpagebreak
\SectionNumber=1\tNumber=1
\flushpar{\bf \number\SectionNumber.\quad  Rees algebras and rational normal scrolls.}
\medskip
Let $k$ be a field, $R$  the standard graded polynomial ring $R=k[x,y]$, and 
  $I$  a height two ideal of $R$  which is minimally generated by $m$ homogeneous forms of degree $d$. The Hilbert-Burch Theorem guarantees that $I$ is generated by the maximal order minors of an $m\times (m-1)$ matrix $\f$, with homogeneous entries. In this paper,  the ideal $I$ is ``almost linearly presented'' in the sense that
the entries of one column of $\f$ have degree $n$; all of the other entries of $\f$ are linear. So, $d=n+m-2$ and the resolution of $I$ looks like  
$$0\to \matrix R(-d-1)^{m-2}\\\p\\R(-d-n)\endmatrix @> \f>> R(-d)^m @>\bmatrix \delta_1&\dots&\delta_m\endbmatrix >>  I\to 0.\tag\tnum{Goal} $$ 
The Rees algebra of $I$ is equal to $\Cal R(I)=R[It]$. Let $S$ be the polynomial ring 
$S=R[T_1,\dots,T_m]$ in $m$ indeterminates over $R$ and let $\Cal A$ be  the kernel of the $R$-algebra homomorphism
 $\Pi\:S\to \Cal R(I)$ which sends $T_i$ to $\delta_it$. 
In this section we identify an $S$-ideal $H$ so that  $A=S/H$ is a  normal domain  and $\Cal AA$ is a height one ideal of $A$. (The ring $A$ is   the coordinate ring of a rational normal scroll.) We also identify an explicit divisorial ideal $K^{(n)}$ of $A$ which is generated by monomials  and an explicit element $g$ of $S$. The main result of the present section is  Theorem \tref{main1} where we prove that 
the ideals $y^n\Cal AA$ and $gK^{(n)}$ of $A$ are equal.

We identify an explicit generating set for $\Cal A$ in Theorem \tref{A=L}.  An explicit minimal generating set for the ideal $K^{(n)}$ may be found in Theorem \tref{Kupn}.

Assume $n\ge 2$. Let $\f'$ denote the restriction of $\f$ to $R(-d-1)^{m-2}$. We call $\f'$ the linear part of $\f$ and we see that the image of $\f'$ is the 
$R$-submodule of $\operatorname{syz}_1^R(I)$ which is generated by the component of degree $d+1$, where $\operatorname{syz}_1^R(I)$ is the first syzygy module of the $R$-module $I$.  In other words,
$$\f': R(-d-1)^{m-2}\to [\operatorname{syz}_1^R(I)]_{d+1}R$$ is an isomorphism. The row space of $\f'$, $\operatorname{RowSp}(\f')$, is the $R$-module generated by the rows of $\f'$. Observe that the minimal number of generators of the $R$-module $\operatorname{RowSp}(\f')$, denoted $\mu(\operatorname{RowSp}(\f'))$, depends only on $I$ and not on 
the set of generators $\{\delta_i\}$ for $I$ or the presenting matrix $\f$ of $I$. 
\definition{Definition \tnum{row}} Given the resolution (\tref{Goal}), with $n\ge 2$, let $\rho=\rho(I)$ be the parameter $$\rho=\mu(\operatorname{RowSp}(\f'))-m+2.$$\enddefinition

The hypothesis that $I$ has height two ensures that $m-1\le \mu(\operatorname{RowSp}(\f'))$; and therefore, 
$$1\le \rho\le 2.$$ We have introduced the parameter $\rho$ at the present time for expository reasons; however, ultimately, $\rho$ plays a significant role in our study. For example, the value $\rho$ determines whether the special fiber ring of $I$ is Cohen-Macaulay (see Theorem \tref{T1}). It also determines many analytic properties of the powers of the ideal $I$, see Sections 4 and 6. 

Start with some minimal resolution for $I$:
$$0\to F_{1,1}\p F_{1,2} \to F_0 \to I,$$ with $F_0\cong R(-d)^m$, $F_{1,1}\cong R(-d-1)^{m-2}$ and $F_{1,2}\cong R(-d-n)$. We prove in Proposition \tref{P28.1} 
that there exists 
a partition $\pmb\sigma$ of $m-2$ into $\rho$    pieces  and there exist bases for $F_0$ and $F_{1,1}$ such that  the linear part of $\f$ is equal to the $m\times(m-2)$ matrix 
$$\f'=\cases \bmatrix D_{\sigma_1}&0\\0&D_{\sigma_2}\endbmatrix,&\text{if $\rho=2$,}\\
\bmatrix D_{\sigma_1}\\0\endbmatrix,&\text{if $\rho=1$,}\endcases\tag\tnum{phi'}
$$where 
$D_a$ is the $(a+1)\times a$ matrix
$$D_a=\bmatrix 
x&0&0&0\\
-y&x&0&0\\
0&-y&\ddots&0\\
0&0&\ddots&x\\
0&0&0&-y
\endbmatrix,$$ 
and  
$$\cases\text{$\pmb\sigma=(\sigma_1,\sigma_2)$, with $\sigma_1\ge\sigma_2\ge 1$, and $\sigma_1+\sigma_2=m-2$,}&\text{if $\rho=2$, or}\\
\text{$\pmb\sigma=\sigma_1=m-2$,}&\text{if $\rho=1$.}\endcases$$ 
 We give the variables $T_j$ of $S$ alternate names.
Let 
$$\cases T_{1,j}=T_j,&\text{if $1\le j\le \sigma_1+1$, and}\\
T_{2,j}=T_{\sigma_1+1+j},&\text{if $\rho=2$ and $1\le j\le \sigma_2+1$.}\endcases\tag\tnum{Alt}$$
\definition{Definition \tnum{CN}}Adopt the data of (\tref{Goal}) with $\f=\bmatrix \f'&\f''\endbmatrix$, where $\f'$ is given in (\tref{phi'}) and $\f''$ is an  $m\times 1$ matrix  of homogeneous forms of degree $n$.
Let 
$\psi$ be the
$2\times (m-1)$ matrix  
$$\psi=\cases {\bmatrix \psi_1&\psi_2&\psi_{3}
\endbmatrix,}&\text{if $\rho=2$,}\\
\bmatrix \psi_1&\psi_2\endbmatrix,&\text{if $\rho=1$},\endcases$$where
each $\psi_i$ is a generic scroll matrix:
$$\psi_i=\cases\bmatrix T_{i,1}&T_{i,2}& \dots&T_{i,\sigma_i-1}&T_{i,\sigma_i} \\T_{i,2}&T_{i,3}&\dots&T_{i,\sigma_i}&T_{i,\sigma_i+1} \endbmatrix&\text{if $1\le i\le \rho$}\\
\bmatrix y\\x\endbmatrix&\text{if $i=\rho+1$}\endcases\tag\tnum{psi}$$
Let $H=I_2(\psi)$,  $A=S/H$,
 $\pmb T$ be the matrix $\bmatrix T_1&\dots& T_m\endbmatrix$,   $g\in S$ be the 
product   $\pmb T\f''$, and $K$ be the ideal of $A$ which is generated by the entries in the top row of $\psi$. The ring $A$ is a domain; let   $Q$ be the quotient field of $A$.  The ideal $K$ is a prime ideal of $A$; let $K^{(n)}$ be the $n^{\text{\rm th}}$ symbolic power of $K$.
\enddefinition
  
\bigskip 
Observe  that  the $A$ is a normal domain of dimension four and $K$ is a height one prime ideal of $A$. It is convenient to think of the ring $S$ as bi-graded. 
$$\matrix \text{The variables $\{T_i\}$ have degree $(0,1)$.}\hfill\\ \text{The variables $\{x,y\}$ have degree $(1,0)$.}\hfill\endmatrix\tag\tnum{bi-}$$ 
Notice that $H$ is a homogeneous ideal with respect to this bi-grading   
and thus (\tref{bi-}) induces a grading on $A$.
The last column of  $\pmb T \f$
has the form $$g=\sum\limits_{i=0}^n c_ix^{n-i}y^i\in S, \tag\tnum{g}$$
where  $c_0,\dots,c_n$ are homogeneous elements of $S$ of degree $(0,1)$ and  $g$ is  a homogeneous 
polynomial of degree $(n,1)$. The generators of $\Cal A$ which are not   in $H$  are all described in terms of the polynomials $c_0,\dots,c_n$; see   Definition \tref{1D5.3}.

\remark{Remark} Let $V\subseteq \Bbb P^{m+1}$ be the variety defined by $I_2(\psi)$. We observe that when $\rho=1$, then the defining equations of $V$ do not involve the variable $T_m$. In other words, in this case, $V$  is the cone over a rational normal scroll $V'\subseteq \Bbb P^m$. If $\rho=2$, then $V$ itself is a rational normal scroll.\endremark

\proclaim{Observation \tnum{O5.1}} The ideal $H$ and the polynomial $g$ are contained in $\Cal A$.\endproclaim
\demo{Proof}
The symmetric algebra $\Sym(I)$ is equal to 
$$\frac S{I_1\(\pmb T  \f\)},$$and the homomorphism 
$\Pi\:S\to \Cal R(I)$ factors through the natural quotient map $S\to\Sym(I)$; so, 
$$I_1\(\pmb T \f\)\subseteq \Cal A.$$
In particular, $g=\pmb T \f''$ is in $\Cal A$. Write $\psi=\bmatrix \psi'&\psi''\endbmatrix$, where $\psi'$ consists of the first $m-2$ columns of $\psi$ and $\psi''$ is the final column of $\psi$. Observe that the product $\pmb T \f'$ is also 
 equal to 
$$\bmatrix x&-y\endbmatrix \psi'.\tag\tnum{here}$$ Each entry of the matrix (\tref{here})  is equal to a  $2\times 2$ minor of $\psi$ which involves the 
last column. Let $\delta$ be a $2\times 2$ minor of $\psi'$. Since the entries of the product matrix (\tref{here}) are in $\Cal A$, 
Cramer's rule shows that $(x,y)\delta\subseteq \Cal A$; but the ideal $\Cal A$ is 
prime and $\Cal A\cap R=\{0\}$; so $\delta$ is also in $\Cal A$. \qed\enddemo

Notice that  $\Cal A A$ is a prime ideal of height one in $A$ because the Rees algebra $\Cal R(I)\cong A/\Cal AA$ is a domain of dimension three.  
 
\proclaim{Theorem \tnum{main1}} Retain the data of Definition \tref{CN}. The following statements hold.
\roster
\item"{(a)}"The ideals $y^n\Cal AA$ and $gK^{(n)}$ of $A$ are equal.
\item"{(b)}"The elements $y^n$ and $g$ of $S$ are not in $H$.
\item"{(c)}"The bi-graded $A$-modules $\Cal AA$ and $K^{(n)}(0,-1)$ are isomorphic. {\rm(}The grading is described in 
   {\rm(\tref{bi-}))}.
\endroster\endproclaim

\demo{Proof}Degree considerations show that $y$ is not in $H$. The ideal $H$ is prime, so $y^n$ is also not in $H$. Assertion (c) will follow from (a) and (b) because $y^n$ has bi-degree $(n,0)$ and $g$ has bi-degree $(n,1)$.
Write $\overline{\phantom{x}}$ to mean image in $A$. We prove (a) by showing that  $\overline{\Cal A}=({\bar g}/{\bar y^n})K^{(n)}$, where the fraction is taken in  $Q$.

We first claim that 
$$(\bar y^i)\:\!\!_QK^{(i)}=(\bar x,\bar y)^i\tag\tnum{1.12}$$ for all $i\ge 1$. Since $\bar y\in K$, one has
$(\bar y^i)\link\!_QK^{(i)}\subseteq (\bar y^i)\link\!_Q (\bar y^i)=A$. Therefore, $(\bar y^i)\link\!_QK^{(i)}=(\bar y^i)\link\!_AK^{(i)}$. The determinantal relations of $A$ give $(\bar x,\bar y)K\subseteq (\bar y)$. Raise each side to the $i^{\text{th}}$ power to obtain $(\bar x,\bar y)^{i}K^{i}\subseteq (\bar y^i)$. Now, it is not difficult to see that
 $(\bar x,\bar y)^{(i)}K^{(i)}\subseteq (\bar y^i)$. (Keep in mind that one need only verify the inclusion locally at all associated prime ideals of $(\bar y^i)$; every associated prime of $(\bar y^i)$ has height one because $A$ is a Cohen-Macaulay domain; and every localization of $A$ at a height one prime ideal is a DVR because $A$ is a normal domain.)  
Therefore, $(\bar x, \bar y)^{(i)}\subseteq (\bar y^i)\link\!_AK^{(i)}$. This inclusion is an equality because $(\bar x,\bar y)$ is a prime ideal not containing $K$. 
Thus, we have shown that $(\bar y^i)\:\!\!_QK^{(i)}=(\bar x,\bar y)^{(i)}$.
Temporarily giving $x$ and $y$ degree $1$ and the variables $T_i$ degree $0$, we see that $\operatorname{gr}_{(\bar x,\bar y)}(A)\iso A$, which is a domain. Therefore, $(\bar x, \bar y)^{(i)}=(\bar x, \bar y)^i$. This completes the proof of (\tref{1.12}). 

We have $\bar g\in(\bar x, \bar y)^n=(\bar y^n)\link K^{(n)}$, where the last equality holds by (\tref{1.12}). Thus, $\bar gK^{(n)}\subseteq \bar y^n A$. Define $L$ to be the ideal   $(\bar g/\bar y^n)K^{(n)}$ of $A$.   At this point, we see that the ideal $L$ is either zero or divisorial. 

To show that $L$ is not zero and to establish  
the equality $\overline{\Cal A}=L$, it suffices to prove that $\overline{\Cal A}\subseteq L$, because $\overline{\Cal A} $ is a height one  prime ideal of $A$. Notice that $\bar g\in L$ as $\bar y\in K$. For every $w\in (x,y)R$, one has $I_w=R_w$. Therefore, $\Sym(I)_w=R[It]_w$ and we obtain $(\bar g)_w=\overline{\Cal A}_w$. It follows that $\bar g\neq 0$ (which completes the proof of (b))  and   $\overline{\Cal A}_w\subseteq L_w$. 

To complete the proof of the inclusion  $\overline{\Cal A}\subseteq L$, it suffices to show that some $w\in (x,y)R$ is a non zerodivisor modulo $L$, equivalently, $(\bar x,\bar y)\overline{R}$ is not contained in 
the union of all associated primes of $L$, or yet equivalently, $(\bar x,\bar y)=(\bar x,\bar y)A$ is not contained in any associated  prime of $L$. But this simply means that $L\not \subseteq (\bar x, \bar y)$, because the ideal $L$ is divisorial and $(\bar x, \bar y)\neq 0$ is prime. 

Finally, to show that $(\bar g/\bar y^n)K^{(n)}=L\not \subseteq (\bar x, \bar y)$, we compute $(\bar x, \bar y)\link\!_QK^{(n)}$ and verify that this fractional ideal does not contain $\bar g/\bar y^n$. Using (\tref{1.12}) twice, we deduce

$$\split &(\bar x, \bar y)\link\!_Q K^{(n)}= \((\bar y)\link\!_Q K\)\link\!_Q K^{(n)}=(\bar y)\link\!_Q KK^{(n)}=\bar y^{-n}[(\bar y^{n+1})\link\!_Q KK^{(n)}]\\
&{}=\bar y^{-n}[(\bar y^{n+1})\link\!_Q K^{n+1}]= \bar y^{-n}(\bar x, \bar y)^{n+1}.\endsplit $$
This fractional ideal cannot contain $\bar g/\bar y^n$; for otherwise, $\bar g\in (\bar x, \bar y)^{n+1}$, which is impossible because $\bar g$ is a non-zero homogeneous element of degree $(n,1)$.
\qed\enddemo

\bigskip
\bigskip
\bigpagebreak
\SectionNumber=2\tNumber=1
\flushpar{\bf \number\SectionNumber.\quad  Matrices with linear entries.}
\medskip

Let $\f$ be the matrix of (\tref{Goal}). In Proposition \tref{P28.1} we prove that there exist row and column operations on $\f$ which transform the linear part of $\f$ into a matrix of the form described in (\tref{phi'}). Recall that $R$ is  the polynomial ring $k[x,y]$ over the field $k$.  For each non-negative integer ${\sigma}$, let 
$D({\sigma})$ be the $({\sigma}+1)\times {\sigma}$ matrix with 
$$D({\sigma})_{i,j}=\cases 
x,&\text{if $i=j$ and $1\le j\le {\sigma}$,}\\
y,&\text{if $i=j+1$ and $1\le j\le {\sigma}$, and}\\
0&\text{otherwise.}\endcases$$ We see that $D(0)$ is invisible, 
$$D(1)=\bmatrix x\\y\endbmatrix,\and D(2)=\bmatrix x&0\\y&x\\0&y\endbmatrix.$$The matrix $D_{{\sigma}}$ of Section 1 is the same as the matrix $D({\sigma})$ of the present section, with $y$ replaced by $-y$. One may use elementary row and column operations to transform either one of these matrices into the other one.

\proclaim{Proposition \tnum{P28.1}}Let $M$ be an $m\times (m-2)$ matrix whose entries are homogeneous linear forms from $R$. Suppose that 
there exists a column vector $\f''$ in $R^m$ of homogeneous forms of the same degree,  such that the ideal of maximal minors of $\bmatrix M&\f''\endbmatrix$ is an ideal of height two in $R$.  Then there exist matrices  $U\in \GL_m(k)$ and $V\in \GL_{m-2}(k)$ and non-negative integers ${\tau}\le {\sigma}$,   with ${\sigma}+\tau=m-2$, such that $UMV$ is equal to
$$\bmatrix D({\sigma})&0\\0&D({\tau})\endbmatrix.\tag\tnum{dab}$$ 
\endproclaim
\remark{Remark}   If ${\tau}=0$ and ${\sigma}=m-2$, then the matrix  of (\tref{dab}) is
$$\bmatrix D(m-2)\\0\endbmatrix,$$ where $0$ represents a $1\times m-2$ matrix of zeros. Observe that $\mu(\operatorname{RowSp}(M))=m-1$ and the parameter $\rho$ of (\tref{row}) is $1$. In the language of (\tref{phi'}), this is the situation in which  the partition $\pmb \sigma$ of $m-2$ consists of  $1$ piece $\pmb \sigma=(m-2)$.\endremark

\bigskip The proposition follows from the next two lemmas. Lemma \tref{'a+b} shows that the hypotheses of Proposition \tref{P28.1} imply the assumption of Lemma \tref{2.7b}.

\proclaim{Lemma \tnum{'a+b}}Let $\f$ be an $m\times (m-1)$ matrix with entries from some commutative  ring. Suppose that there are positive  integers $p$ and $q$ with $p+q=m$ and
$$\f=\bmatrix Z&Y\\X&W\endbmatrix,$$ where $Z$ is an $p\times q$ matrix of zeros and $Y$, $X$, and $W$ are matrices. Then the ideal $I_{m-1}(\f)$ is contained in the principal ideal 
$(\det X)$. \endproclaim

\demo{Proof} One sees this by expanding any maximal minor along the submatrix consisting of its first $q$ columns. \qed\enddemo

\proclaim{Lemma \tnum{2.7b}} Let $M$ be an $m\times (m-2)$ matrix whose entries are homogeneous linear forms from $R$. Assume that for every  $$\matrix \text{$X\in \GL_m(k)$  and $Y\in \GL_{m-2}(k)$, the product matrix   $XMY$ does not}\hfill\\\text{contain a $p\times q$ submatrix of zeros for any pair of positive integers $(p,q)$}\hfill\\\text{with $p+q=m$.}\hfill\endmatrix\tag\tnum{cond}$$
Then there exist matrices  $U\in \GL_m(k)$ and $V\in \GL_{m-2}(k)$ and non-negative integers ${\tau}\le {\sigma}$,   with ${\sigma}+\tau=m-2$, such that $UMV$ is given in 
 {\rm (\tref{dab})}. 
\endproclaim

\demo{Proof}
The proof   is by induction on $m$. The assertion is obvious when $m=3$. Henceforth, $4\le m$. Let $\overline{M}$ be the image of $M$ in the ring $R/(y)$. We see that $\overline{M}=xM'$ for some $m\times (m-2)$ matrix $M'$ with entries in $k$. There exist invertible matrices 
$X$ and $Y$ with entries in $k$ so that $$XM'Y=\bmatrix M''\\ 0_{2\times (m-2)}\endbmatrix,$$for some matrix $M''$.  Therefore every entry in the bottom two rows of $XMY$ is in the ideal $(y)$. Some entry of the bottom two rows of $XMY$ is not zero by (\tref{cond}). Thus, further row and column operations yield a matrix of the form
 $$
\bmatrix M_1&  M_2\\
 0 &y\endbmatrix. $$
The $(m-1)\times (m-3)$ matrix  $M_1$ satisfies (\tref{cond})  because if there exist $X_1\in \GL_{m-1}(k)$ and $Y_1\in \GL_{m-3}(k)$ so that $X_1M_1Y_1$ contains an $p_1\times q_1$ zero submatrix, then there exist invertible matrices $X$ and $Y$ so that $XMY$ contains an $(p_1+1)\times q_1$ zero submatrix. 
By induction $M$ may be transformed into $$
 \bmatrix D({\sigma})& 0 &C_1\\
0&D({\tau})&C_2\\
0&0&y\endbmatrix. \tag\tnum{form}$$  for two non-negative integers ${\tau}\le {\sigma}$ with ${\sigma}+{\tau}=m-3$, where $C_1$ and $C_2$ are column vectors. Use column operations to remove all $x$'s from $C_1$ and $C_2$, except possibly in the bottom row. Use row operations to remove all $y$'s from $C_1$ and $C_2$. Thus, $M$ may be transformed into a matrix of the form (\tref{form}) with $$C_i=\bmatrix 0\\\vdots\\0\\c_ix\endbmatrix\tag\tnum{ci}$$  for some $c_i\in k$. At least one of the constants $c_1$ or $c_2$ must be non-zero. If $c_1$ is not zero, then pre-multiply and post-multiply by $$U=\bmatrix c_1^{-1}I_{\sigma+1}&0&0\\0&I_{\tau+1}&0\\0&0&1\endbmatrix \and V=\bmatrix c_1I_{\sigma}&0&0\\0&I_{\tau}&0\\0&0&1\endbmatrix, $$ respectively, to transform $c_1$ into $1$. The constant $c_2$ may be treated in a similar manner. Thus, $M$ may be transformed into a matrix of the form (\tref{form}), 
where the columns $C_1$ and $C_2$ are described in (\tref{ci}), and one of the following three cases occurs:
$$\cases c_1=1, c_2=0&\text{case 1}\\c_1=0, c_2=1&\text{case 2}\\c_1=c_2=1&\text{case 3}.\endcases$$
The third case may be transformed into the second case using
$$U=\bmatrix I_{{\sigma}-{\tau}}&0&0&0\\0&I_{{\tau}+1}&-I_{{\tau}+1}&0\\0&0&I_{{\tau}+1}&0\\0&0&0&1\endbmatrix\and V=\bmatrix I_{{\sigma}-{\tau}}&0&0&0\\0&I_{{\tau}}&I_{{\tau}}&0\\0&0&I_{{\tau}}&0\\0&0&0&1\endbmatrix.$$ In the second case, (\tref{form}) is readily seen to be $$\bmatrix D({\sigma})&0\\0&D(\tau+1)\endbmatrix,\tag\tnum{**}$$  and in the first case, one  may rearrange the rows and columns of (\tref{form}) to obtain 
$$\bmatrix D({\sigma}+1)&0\\0&D({\tau})\endbmatrix.$$ Finally, we notice that if ${\tau}+1>{\sigma}$, then one may rearrange the rows and columns of (\tref{**}) to obtain $$\bmatrix D({\tau}+1)&0\\0&D({\sigma})\endbmatrix. \qed$$ 
\enddemo

Proposition \tref{P28.1} shows that any ideal $I$ as in described in (\tref{Goal}) has a presentation matrix in which the linear part is given in (\tref{dab}). One may use elementary row and column operations to transform a matrix given in (\tref{dab}) to a matrix given in   (\tref{phi'}), or vice versa.

The next result was obtained during a conversation with  
David Eisenbud.
\proclaim{Corollary \tnum{C2.4}}If $I$ is a height two ideal in $R=k[x,y]$, then the resolution of $I$ is given in {\rm(\tref{Goal})} if and only if there exists non-negative integers ${\sigma}$ and ${\tau}$, with ${\sigma}+{\tau}= m-2$, and relatively prime homogeneous forms $F_1$ and $F_2$ in $R$, with $\deg F_1=n+{\sigma}$ and $\deg F_2=n+{\tau}$ such that
$$I=(x,y)^{{\tau}}F_1+(x,y)^{{\sigma}}F_2.\tag\tnum{spot}$$\endproclaim
\demo{Proof} Start with the data ${\sigma}$, ${\tau}$, $F_1$, and $F_2$. Write $$F_1=\sum\limits_{i=0}^{{\sigma}}\a_ix^{{\sigma}-i}y^i\and F_2=\sum\limits_{i=0}^{{\tau}}\b_ix^{{\tau}-i}y^i,$$ for homogenous forms $\a_i$ and $\b_i$ of degree $n$. Let
$$\pmb \a=\bmatrix \a_{{\sigma}}\\\vdots\\\a_{0}\endbmatrix\and \pmb \b=\bmatrix \b_{{\tau}}\\\vdots\\\b_{0}\endbmatrix.$$ Observe that
$$\det \bmatrix D_{\sigma}&\pmb \a\endbmatrix =F_1,\quad \det \bmatrix D_{{\tau}}&\pmb \b\endbmatrix=F_2,$$ and the ideal generated by the  maximal order minors of the matrix
$$\bmatrix  D_{{\sigma}}&0&\pmb \a\\
0& D_{{\tau}}&\pmb \b
\endbmatrix$$ is equal to $I$. 

The converse follows from Proposition \tref{P28.1}. \qed  \enddemo
\remark{Remark} The ideal $I$ of (\tref{spot}) is a truncation of a complete intersection: $I=(F_1,F_2)_{\ge d}$. That is, $I$ is generated by all homogeneous elements of the complete intersection $(F_1,F_2)$ of degree at least $d$, where, as always, $d=n+m-2$.\endremark

We remark that H\`a T\`ai has previously studied the Rees algebra of a truncation of an ideal. Let $J$ be the defining ideal of a finite set $\Bbb X$ of points in $\Bbb P^2$ and let $\a$ represent the minimal degree of a generator of $J$. It is shown in \cite{\rref{Ha}} that the Rees algebra $\Cal R(J_{\a+1})$ is Cohen-Macaulay when $\Bbb X$ is a general set of points, and,  for an arbitrary set of points, $\Cal R(J_{t})$ is Cohen-Macaulay for all sufficiently large $t$. In each case the degrees of the generators of the defining ideal of the Rees algebra are given. Our Rees algebras are never Cohen-Macaulay.

\example{Example \tnum{E2.5}} If $F_1=y^{n+{\sigma}}$ and $F_2=x^{n+{\tau}}$, then $I$ is the monomial ideal
$$(y^d,xy^{d-1},\dots,x^{\tau}y^{d-{\tau}})+(x^{d-{\sigma}}y^{\sigma},\dots,x^{d-1}y,x^d).$$\endexample

The following proof was prompted to us by a question of Craig Huneke.
\proclaim{Corollary \tnum{C2.6}} Adopt the notation of Corollary \tref{C2.4} with ${\tau}\le {\sigma}$, and write $d=n+{\sigma}+{\tau}$. For every $s \geq 1$ one has
$${\reg} I^s =  {\max} \{sd, sd -(s-1){\tau}+n-1\}\, .$$
\endproclaim
\demo{Proof} Write $\maxm =(x,y)$.  Notice that the regularity of a homogeneous $\maxm$-primary ideal is the smallest power of $\maxm$ contained in it.
Notice that $I^s$ is generated by forms of degree $sd$ and
$$I^s=\sum\limits_{i=0}^{s} \maxm^{sd-\deg(F_1^iF_2^{s-i})}F_1^iF_2^{s-i} = (F_1,F_2)^s \cap \maxm^{sd} \, .$$
Hence $\maxm^t \subseteq I^{s}$ if and only if  $\maxm^t \subseteq (F_1,F_2)^s$ and $t \geq sd$. In other words, $$\reg I^s = \max \{sd, \reg (F_1,F_2)^s\}\, .$$ Finally, $F_1, F_2$ are a regular sequence of forms of degrees $n+{\tau} \le n+{\sigma}$. Hence $(F_1,F_2)^s$ is presented by the $s+1$ by $s$ matrix
$$\bmatrix
F_2& & & & & \\
-F_1&F_2 &   & & &  \\
&-F_1& \cdot & & &   \\
& &\cdot &  \cdot & &  \\
& & &\cdot & \cdot&  \\
& & & & \cdot& F_2 \\
& & & & & -F_1
\endbmatrix$$
From this minimal homogeneous resolution one sees that
$$ \reg(F_1,F_2)^s= s(n+{\sigma})+n+{\tau}-1=sd-(s-1){\tau}+n-1 . \qed $$ \enddemo

It is shown in \cite{\rref{CHT},\rref{Ko},\rref{TW}} that the  regularity 
of the $s^{\text{\rm th}}$ power of any homogeneous ideal is a linear function of $s$ for all $s\gg 0$. Indeed, in our notation, the aforementioned papers guarantee that $\operatorname{reg}(I^s)=sd+e$ for some non-negative integer $e$. The integer $e$ has been determined in \cite{\rref{EH}}. From Corollary \tref{C2.6}, we  read the exact value of $e$ and   the exact values of $s$ for which the above equation  holds. The answers   depend  on the value of $\rho$.
In Section 6 we resolve each power $I^s$; thereby confirming the present calculations, see especially Corollary \tref{reg}.

\proclaim{Corollary \tnum{C2.7}} Let $I$ be the ideal of Definition \tref{CN} and $s$ be a positive integer.
\flushpar{\bf (1)} If $\rho=1$, then   $\reg I^s= sd+n-1$ for all $s\ge 1$.
\medskip \flushpar{\bf (2)} If $\rho=2$, then  $\reg I^s= sd$ if and only if  $\frac{n-1}{\sigma_2}+1\le s$.\endproclaim

\demo{Proof}If $\rho=1$, then the parameter ${\tau}$ of Corollary \tref{C2.6} is equal to zero and $\max\{sd, sd-(s-1){\tau}+n-1\}$ is equal to $sd+n-1$ for all $s\ge 1$. If $\rho=2$, then the parameter ${\tau}$ of Corollary \tref{C2.6} is equal to $\sigma_2$ and
$\reg I^s=sd$ if and only if $sd\ge sd-(s-1)\sigma_2+n-1$.
\qed\enddemo

In Section 3 we calculate an explicit generating set for the ideal $\Cal A$ which defines the Rees algebra $\Cal R(I)$ for $I$ given in (\tref{Goal}). An alternate approach to this problem is suggested by Corollary \tref{C2.4}.
\remark{Remark \tnum{E53.13}}
Let $I$ be the ideal $\maxm^{\tau}F_1+\maxm^{\sigma}F_2$ of $R=k[x,y]$, as described in Corollary \tref{C2.4}, where $\maxm$ is the maximal homogeneous ideal $(x,y)$ of $R$ and $F_1,F_2$ is a regular sequence of homogeneous forms in $R$ with $\deg F_1 =n+\sigma$ and $\deg F_2=n+\tau$. The ideal $I$ is contained in the complete intersection ideal $(F_1,F_2)$ and the Rees algebra of $(F_1,F_2)$ is well understood. Let $D=R[u,v]$ and  
$R[t]$ be   polynomial rings 
 and $M_2\: D\to R[t]$ be the $R$-algebra homomorphism with 
$M_2(u)=F_2t$ and $M_2(v)=-F_1t$. 
The image of $M_2$ is the Rees algebra $\Cal R((F_1,F_2))=R[F_1t,F_2t]$.
We have a short exact sequence
$$0\to (\Phi)\to D @> M_2 >> R[F_1t,F_2t]\to 0,$$ where $\Phi=F_1u+F_2v$. 

Define $C=R[R_{\sigma}u,R_{\tau}v]\subseteq D$ and
$M_1\: S\to C$, with $M_1(T_{1,j})=ux^{j-1}y^{\sigma+1-j}$ and 
$$\cases M_1(T_{2,j})=vx^{j-1}y^{\tau+1-j}&\text{if $\rho=2$}\\
M_1(T_m)\phantom{_{,}}=v&\text{if $\rho=1$.}\endcases$$ Observe that $M_1$ induces an isomorphism $A\to C$, which we also call $M_1$. 
The image of the restriction of $M_2$ to $C$ is equal to the Rees algebra $\Cal R(I)=R[It]$. The ideal $\Cal AA$ is equal to 
$$\ker \(A@> M_1>> C@> M_2>> R[It]\).$$
The commutative diagram
$$\CD 0@>>> (\Phi)@>>> D@> M_2>> R[F_1t,F_2t]@>>> 0\\
@. @. @A \incl AA @A \incl AA\\
0@>>>(\Phi)D\cap C @>>> C@> M_2|_C >> R[It]@>>>0\endCD$$has exact rows.
Thus,
$$ \Cal A A@> M_1>> (\Phi)D\cap C$$is an isomorphism.
It is easy to see that 
$$(\Phi)D\cap C=\sum\limits_{i,j}(x,y)^{i\sigma+j\tau-n}\Phi u^iv^j,$$where the sum is taken over all non-negative integers $i,j$ with $i\sigma+j\tau-n\ge 0$. The problem of determining a generating set for the defining ideal of $R[It]$ is equivalent to the problem lifting the generators of 
$(\Phi)D\cap C$ to $A$. (For example, in the present notation, one can check that $M_1^{-1}(\Phi)$ is equal to the image of   $g$ in $A$   because
$$g=\sum_{i=0}^{\sigma}T_{1,\sigma+1-i}\a_i+\sum_{i=0}^{\tau}T_{2,\tau+1-i}\b_i,$$ see the proof of Corollary \tref{C2.4}.)  This is a non-trivial calculation, similar in difficulty to the calculation of Section 3. Our approach in Section 3 is to  give $\Cal AA$ the structure of a divisor on a scroll. The advantage of the divisorial approach is that we have an explicit isomorphism between $\Cal AA$  and a well understood monomial ideal.  We determine a minimal generating set of $\Cal AA$, the degrees of the minimal generators, and even a resolution  of $R[It]$.   
\endremark 

\bigskip
\bigskip
\bigpagebreak
\SectionNumber=3\tNumber=1
\flushpar{\bf \number\SectionNumber.\quad  Explicit  generators   for the defining ideal of the Rees algebra.}
\medskip

The main result of this section is Theorem \tref{A=L} where we identify an explicit generating set for the defining ideal $\Cal A$ of the Rees algebra $\Cal R(I)$. 
Adopt the data of Definition \tref{CN} with
$$\ell=\rho+1,\quad\sigma_{\ell}=1,\quad y=T_{\ell,1},\and x=T_{\ell,2}.\tag\tnum{3.1}$$ In this notation, the matrix $\psi_{\ell}$ of (\tref{psi}) is 
$$\psi_{\ell}=\bmatrix T_{\ell,1}\\T_{
\ell,2}\endbmatrix.$$According to Theorem \tref{main1}, we need to identify generators for the ideal $\Cal L$ in $S$ with $y^n\Cal LA=gK^{(n)}$. The following minimal generating set for $K^{(n)}$ is calculated in   \cite{\rref{KPU-d}, Prop\.~1.20}. 

\proclaim{Theorem \tnum{Kupn}} 
A $k$-tuple $\pmb a=(a_1,\dots,a_k)$ of non-negative integers is {\it eligible} if $0\le k\le \rho$ and $\sum\limits_{u=1}^k a_u\sigma_u<n$.   If $\pmb a$ is an eligible $k$-tuple, then 
$f(\pmb a)$ and $r(\pmb a)$ are defined by{\rm :}
$$\sum\limits_{u=1}^ka_u\sigma_u+f(\pmb a)\sigma_{k+1}<n\le \sum\limits_{u=1}^ka_u\sigma_u+(f(\pmb a)+1) \sigma_{k+1}$$and
$$r(\pmb a)= \sum\limits_{u=1}^ka_u\sigma_u+(f(\pmb a)+1) \sigma_{k+1}-n+1.$$
The ideal $K^{(n)}$ of $A$ is equal to
$$K^{(n)}=(\{T^{\pmb a}T_{k+1,1}^{f(\pmb a)}T_{k+1,j}\mid \text{$\pmb a$ is an eligible $k$-tuple and $1\le j\le r(\pmb a)$}\})A,$$
where $T^{\pmb a}=\prod\limits_{u=1}^kT_{u,1}^{a_u}$.  
\endproclaim

\remark{Remark} The empty tuple, $\emptyset$, is always eligible, and we have
$$\tsize f(\emptyset)=\lceil\frac n{\sigma_1}\rceil -1,\quad r(\emptyset)=\sigma_1\lceil\frac n{\sigma_1}\rceil -n+1, \and T^{\emptyset}=1.$$\endremark

\definition{Definition \tnum{-1D6.1}} Recall the polynomials $c_0,\dots,c_n$ of (\tref{g}). 
\smallskip\flushpar{(a)}  
 For  integers $a$ and $b$ with $a+b\le n$ and $0\le a$, define the polynomial 
 $\Delta_{a,b}$ to be 
$$\cases \sum\limits_{k=0}^bc_{a+k}x^{b-k}y^k=c_ax^b+c_{a+1}x^{b-1}y
+\dots+ c_{a+b}y^b,&\text{if $0\le b$},\\ 0,&\text{if $b<0$.}\endcases$$
In particular $g=\Delta_{0,n}$. Furthermore, $\Delta_{a,b}$ is a homogeneous element of $S$ of degree 
$(b,1)$.
\smallskip\flushpar{(b)} If $0\le a\le n$, then write $\Delta_a$ to mean $\Delta_{a,n-a}$. So
$$\Delta_a=c_ax^{n-a}+c_{a+1}x^{n-a-1}y+\dots+c_ny^{n-a},$$ and
$\Delta_a$ is a homogeneous element of $S$ of degree $(n-a,1)$.
\smallskip\flushpar{(c)} For each $4$-tuple of non-negative indices $(i,a,b,\g)$ with $$1\le i\le 2,\quad b+1\le \g\le \sigma_i+1, \and a+b\le n,$$ define 
$$\pi_{i,a,b,\g}=\sum_{k=0}^{b}c_{a+k}T_{i,\g-k}=c_aT_{i,\g}+c_{a+1}T_{i,\g-1}+\dots+c_{a+b}T_{i,\g-b}.$$
 This element of 
$S$ has bi-degree $(0,2)$ if $i\le\rho$ and bi-degree $(1,1)$ if $i=\ell$. 

\smallskip\flushpar{(d)} If $1
\le i\le 2$   and $0\le a\le n-\sigma_i+1$, then let $\pi_{i,a}$ mean $\pi_{i,a,\sigma_i-1,\sigma_i+1}$; so $\pi_{i,a}$ is equal to
$$ \sum\limits_{k=0}^{\sigma_i-1}c_{a+k}T_{i,\sigma_i+1-k}=c_aT_{i,\sigma_i+1}+c_{a+1}T_{i,\sigma_i}+\dots+c_{a+\sigma_i-1}T_{i,2}.
$$
\smallskip\flushpar{(e)} If $(i,s,j)$ are non-negative integers with $1\le i\le 2$, $s\le n$, and 
$1\le j\le \sigma_i+1-s$, then let $\pi_{i,s,j}'$ mean $\pi_{i,n-s,s,s+j}$.
\enddefinition

\remark{Remarks \tnum{1R6.8}}
{(a)}
 Reverse the order of summation in the polynomial  $\pi_{i,s,j}'$ 
to write
$$\pi_{i,s,j}'=\sum_{k=0}^{s}c_{n-k}T_{i,j+k}=
c_nT_{i,j}+c_{n-1}T_{i,j+1}
+\dots+c_{n-s} T_{i,j+s}.$$
\smallskip\flushpar{(b)} If the non-negative integers $a,b,\g$ satisfy
$a+b\le n$ and $1\le \g\le b$, then 
$$\Delta_{a,b}=x^{b-\g+1}\Delta_{a,\g-1}+y^{\g}\Delta_{a+\g,b-\g}.$$
The polynomial $\Delta_{a,b}$ of $S$ is homogeneous in $x$ and $y$ of degree $b$; 
hence, every 
term in $\Delta_{a,b}$  is divisible by either $y^{\g}$ or $x^{b-\g+1}$.
The formula records the fact that we have already chosen names for the coefficients of $\Delta_{a,b}$ in
$(y^{\g},x^{b-\g+1})$.
  At any rate,
the left hand side is
$$\(c_ax^b+\dots+c_{a+\g-1}x^{b-\g+1}y^{\g-1}\)+\(c_{a+\g}x^{b-\g}y^{\g}+\dots+c_{a+b}y^b\)$$$$=x^{b-\g+1}\(c_ax
^{\g-1}+\dots+c_{a+\g-1}y^{\g-1}\)+y^{\g}\(c_{a+\g}x^{b-\g}+\dots+c_{a+b}y^{b-\g}\),$$which is the right hand side.
\smallskip\flushpar{(c)} 
If $N$ is negative, then the sum $\sum\limits_{a+b=N}$ is zero; if $N$ is a non-negative integer then the sum $\sum\limits_{a+b=N}$ is taken over all pairs of non-negative integers $(a,b)$, with $a+b=N$. 
\smallskip\flushpar{(d)} 
We calculate in $S$. If $s_1$ and $s_2$ are elements of $S$, we write $s_1\equiv s_2$ to mean that $s_1-s_2\in H$.
\endremark

\definition{Definition \tnum{1D5.3}} For each pair $(\pmb a,j)$, where $\pmb a$ is an eligible tuple and $1\le j\le r(\pmb a)$, we define a polynomial $G_{(\pmb a,j)}$ in $S$.

\smallskip\flushpar{(a)} If $1\le j\le r(\emptyset)$, then let
$$G_{(\emptyset,j)}=f_j=T_{1,j+\sigma_1+1-r(\emptyset)}\sum\limits_{p+q=f(\emptyset)-1}T_{1,1}^pT_{1,\sigma_1+1}^{q}\pi_{1,p\sigma_1} 
+T_{1,1}^{f(\emptyset)}\pi'_{1,\sigma_1+1-r(\emptyset),j}.$$
\smallskip\flushpar{(b)} If $(a_1)$ is an eligible $1$-tuple, and $1\le j\le r(a_1)$, then let
$$G_{((a_1),j)}=g_{a_1,j}=\left\{\matrix 
T_{2,j+\sigma_2+1-r(a_1)}T_{2,\sigma_2+1}^{f(a_1)}\sum\limits_{p+q=a_1-1}T_{1,1}^pT_{1,\sigma_1+1}^q \pi_{1,p\sigma_1}\hfill\\ \vspace{5pt}
+T_{1,1}^{a_1}T_{2,j+\sigma_2+1-r(a_1)}\sum\limits_{p+q=f(a_1)-1}T_{2,1}^pT_{2,\sigma_2+1}^{q}\pi_{2,a_1\sigma_1+p\sigma_2}\hfill\\\vspace{5pt}
+ T_{1,1}^{a_1}T_{2,1}^{f(a_1)}\pi'_{2,\sigma_2+1-r(a_1),j}.\hfill\endmatrix\right.$$

  \smallskip\flushpar{(c)} If $\pmb a=(a_1,a_2)$ is an eligible $2$-tuple, then $r(\pmb a)=1$. Let
$$G_{(\pmb a,1)}=h_{a_1,a_2}=\left\{\matrix 
x^{n-a_1\sigma_1-a_2\sigma_2}  T_{2,\sigma_2+1}^{a_2}\sum\limits_{p+q=a_1-1}T_{1,1}^pT_{1,\sigma_1+1}^q \pi_{1,p\sigma_1}\hfill\\\vspace{5pt} 
+x^{n-a_1\sigma_1-a_2\sigma_2}T_{1,1}^{a_1}\sum\limits_{p+q=a_2-1}T_{2,1}^pT_{2,\sigma_2+1}^{q}\pi_{2,a_1\sigma_1+p\sigma_2}\hfill 
\\\vspace{5pt}
+ T_{1,1}^{a_1}T_{2,1}^{a_2}\Delta_{a_1\sigma_1+a_2\sigma_2}.\hfill\endmatrix\right.$$

\smallskip\flushpar{(d)} The ideal $\Cal L$ of $S$ is equal to
$$H+\(\{G_{(\pmb a,j)}\mid \text{$\pmb a$ is an eligible tuple and $1\le j\le r(\pmb a)$ }\}\).$$
\enddefinition

 We are now able to state the main result of this section. The ideal $\Cal A$ which defines the Rees algebra $\Cal R(I)$ was introduced in the first paragraph of Section 1.
\proclaim{Theorem \tnum{A=L}} The ideals $\Cal A$ and $\Cal L$ of the ring $S$ are equal.\endproclaim
\demo{Proof} In light of Observation \tref{O5.1} and Theorem \tref{main1}, we need only show that 
the ideals $gK^{(n)}$ and $y^n\Cal LA$ of $A$ are equal.
This calculation is carried out in Lemma \tref{L46.1}.g. \qed \enddemo

\remark{Remarks \tnum{1R6.1}}
\smallskip\flushpar{(a)}  If $\rho=2$, then 
 $$\hskip-11pt \matrix
f_j\hfill&\text{is homogeneous of degree $(0,f(\emptyset)+2)$},\hfill\\
g_{a_1,j}\hfill&\text{is homogeneous of degree $(0,a_1+f(a_1)+2)$, and}\hfill&\\
h_{a_1,a_2}\hfill& \text{is homogeneous of degree $(f(a_1,a_2)+1,a_1+a_2+1)$}.\hfill\endmatrix$$
 If $\rho=1$, then 
 $$\hskip-11pt \matrix
f_j\hfill&\text{is homogeneous of degree $(0,f(\emptyset)+2)$, and}\hfill\\
g_{a_1,j}\hfill&\text{is homogeneous of degree $(f(a_1)+1,a_1+1)$. }\hfill\endmatrix$$
\smallskip\flushpar{(b)} Let $0^s$ be the $s$-tuple $(0,\dots,0)$. Observe that $G_{(0^{\rho},1)}=g$. Indeed, if $\rho=2$, then $h_{0,0}=\Delta_{0}=g$, and if $\rho=1$, then $$\split g_{0,1}&{}=T_{2,2}\sum\limits_{p+q=n-2}T_{2,1}^pT_{2,2}^q\pi_{2,p}+T^{n-1}_{2,1}\pi_{2,n-1,1,2}\\&{}=x^2\sum\limits_{p+q=n-2}y^px^qc_{p}+y^{n-1}(c_{n-1}x+c_ny)=g.\endsplit$$
\endremark

\proclaim{Observation \tnum{1t6.7}}If $a$, $i$, and $j$ are integers with $0\le a$, $1\le i\le\rho$,   $1\le j$, and $ j+a\le \sigma_i+1$ , then
$ x^aT_{i,j}\equiv y^aT_{i,j+a} $.
\endproclaim
\demo{Proof} The ideal $$I_2\bmatrix T_{i,1}&T_{i,2}&\dots&T_{i,\sigma_i-1}&T_{i,\sigma_i}&y\\T_{i,2}&T_{i,3}&\dots&T_{i,\sigma_i}&T_{i,\sigma_i+1}&x\endbmatrix$$ is contained in $H$. A quick induction completes the proof. \qed\enddemo

\proclaim{Observation \tnum{112.22}} Take $1\le i\le \rho$.
\flushpar{\rm(a)} If $0\le a\le n-\sigma_i+1$, then $T_{i,1}x\Delta_{a,\sigma_i-1}\equiv y^{\sigma_i}\pi_{i,a}$.
\smallskip\flushpar{\rm(b)} If $0\le s\le n$ and $1\le j\le \sigma_i+1-s$, then $T_{i,j}\Delta_{n-s}\equiv y^{s}\pi_{i,s,j}'$.
\endproclaim
\demo{Proof}Use Observation \tref{1t6.7} to see the left hand side of (a) is 
 $$\sum_{k=0}^{\sigma_i-1}c_{a+k}\(x^{\sigma_i-k}T_{i,1}\)y^k\equiv \sum_{k=0}^{\sigma_i-1}c_{a+k}\(y^{\sigma_i-k}T_{i,\sigma_i-k+1}\)y^k =y^{\sigma_i}\sum_{k=0}^{\sigma_i-1}c_{a+k}T_{i,\sigma_i-k+1},$$and this is the right hand side of (a).
In a similar manner, we see that the left hand side of (b) is
$$ \align \sum_{k=0}^{s}c_{n-s+k}(x^{s-k}T_{i,j})y^k&{}\equiv \sum_{k=0}^{s}c_{n-s+k}(y^{s-k}T_{i,j+s-k})y^k =
y^{s}\sum_{k=0}^{s}c_{n-s+k}T_{i,j+s-k}\\&{}= 
y^{s}\pi_{i,n-s,s,j+s}=y^{s}\pi_{i,s,j}'. \qed  \endalign$$ 
\enddemo

\proclaim{Lemma \tnum{L46.1}}
\smallskip\flushpar{\rm (a)} If $(a_1,0)$ and $(a_1+1,0)$ are eligible tuples, then
$T_{1,1}h_{a_1,0}\equiv y^{\sigma_1}h_{a_1+1,0}$.
\smallskip\flushpar{\rm (b)} If $(a_1,a_2)$ and $(a_1,a_2+1)$ are eligible tuples, then
$T_{2,1}h_{a_1,a_2}\equiv y^{\sigma_2}h_{a_1,a_2+1}$.
\smallskip\flushpar{\rm (c)} If $(a_1)$ is eligible, $\rho=2$,  and $1\le j\le r(a_1)$, then
$$T_{2,j}h_{a_1,f(a_1)}\equiv y^{\sigma_2+1-r(a_1)}g_{a_1,j}.$$
\smallskip\flushpar{\rm (d)} If $\rho=2$  and $1\le j\le r(\emptyset)$, then
$T_{1,j}h_{f(\emptyset),0}\equiv y^{\sigma_1+1-r(\emptyset)}f_{j}$.
\smallskip\flushpar{\rm (e)} If $\rho=1$  and $(a_1)$ and $(a_1+1)$ are eligible tuples,
then $T_{1,1}g_{a_1,1}\equiv y^{\sigma_1}g_{a_1+1,1}$. 
\smallskip\flushpar{\rm (f)} If $\rho=1$  and $1\le j\le r(\emptyset)$, then $T_{1,j}g_{f(\emptyset),1}\equiv y^{\sigma_1+1-r(\emptyset)}f_{j}$.
\smallskip\flushpar{\rm (g)} The ideals $gK^{(n)}$ and $y^n\Cal LA$ of $A$ are equal.
\endproclaim
\demo{Proof}To prove (a) we recall that
$$T_{1,1}h_{a_1,0}= 
T_{1,1}x^{\sigma_1}x^{n-(a_1+1)\sigma_1}  \sum\limits_{p+q=a_1-1}T_{1,1}^pT_{1,\sigma_1+1}^q \pi_{1,p\sigma_1}
+ T_{1,1}^{a_1+1}\Delta_{a_1\sigma_1}.$$
 The facts 
$$ x^{\sigma_1}T_{1,1}\equiv y^{\sigma_1}T_{1,\sigma_1+1},\tag\tnum{fact1}$$
$$\Delta_{a_1\sigma_1}=x^{n-(a_1+1)\sigma_1+1}\Delta_{a_1\sigma_1,\sigma_1-1}+y^{\sigma_1}\Delta_{(a_1+1)\sigma_1}, \text{ and }\tag\tnum{fact2}$$
$$T_{1,1}x\Delta_{a_1\sigma_1,\sigma_1-1} \equiv y^{\sigma_1}\pi_{1,a_1\sigma_1} \tag\tnum{fact3}$$
may be found  in Observation \tref{1t6.7}, Remark \tref{1R6.8}(b), and Observation \tref{112.22}(a), respectively. Apply (\tref{fact1}) to the first summand of 
$T_{1,1}h_{a_1,0}$ and (\tref{fact2}) and (\tref{fact3}) to the second summand in order to 
 establish (a). 

The same type of methods are used to prove (b). One uses $T_{2,1}x^{\sigma_2}\equiv y^{\sigma_2}T_{2,\sigma_2+1}$ in the first two summands of $T_{2,1}h_{a_1,a_2}$. In the third summand one uses 
$$ \Delta_{a_1\sigma_1+a_2\sigma_2}
=x^{n-a_1\sigma_1-(a_2+1)\sigma_2} x\Delta_{a_1\sigma_1+a_2\sigma_2,\sigma_2-1}+y^{\sigma_2}\Delta_{a_1\sigma_1+(a_2+1)\sigma_2}.$$
Once again, Observation \tref{112.22}(a) yields 
$$T_{2,1}x\Delta_{a_1\sigma_1+a_2\sigma_2,\sigma_2-1}
\equiv y^{\sigma_2} \pi_{2, a_1\sigma_1+a_2\sigma_2}.$$

We prove (c). Notice that
$$n-a_1\sigma_1-f(a_1)\sigma_2=\sigma_2 +1-r(a_1);$$hence
Observations \tref{1t6.7} and \tref{112.22}(b) yield $$\align T_{2,j}x^{n-a_1\sigma_1-f(a_1)\sigma_2}&{}\equiv y^{\sigma_2 +1-r(a_1)}T_{2,j+\sigma_2+1-r(a_1)}\ \text{and}\\T_{2,j}\Delta_{a_1\sigma_1+f(a_1)\sigma_2}&{}\equiv y^{\sigma_2+1-r(a_1)}\pi'_{2,\sigma_2+1-r(a_1),j}.\endalign$$ The proof of (d) is similar. The equality $$n-f(\emptyset)\sigma_1=\sigma_1+1-r(\emptyset)$$ implies
$$T_{1,j}x^{n-f(\emptyset)\sigma_1}\equiv y^{\sigma_1+1-r(\emptyset)}T_{1,j+\sigma_1+1-r(\emptyset)}\ \text{and}\tag\tnum{fact4}$$  $$T_{1,j}\Delta_{f(\emptyset)\sigma_1}\equiv y^{\sigma_1+1-r(\emptyset)}\pi'_{1,\sigma_1+1-r(\emptyset),j}.\tag\tnum{fact5}$$

We now prove (e) and (f). When $\rho=1$, we have $\sigma_2=1$, $T_{2,1}=y$, $T_{2,2}=x$.
For any eligible $1$-tuple $(a)$ one has
 $n-a\sigma_1=f(a)+1$, and $r(a)=1$. We quickly calculate
$$\pi_{2,a\sigma_1+p}=c_{a\sigma_1+p}x,\quad \text{for $0\le p\le f(a)+1$, and}\quad\pi'_{2,1,1}=c_{n-1}x+c_ny.$$ We now have 
$$x\sum\limits_{p+q=f(a)-1}y^px^q \pi_{2,a\sigma_1+p}+y^{f(a)}\pi'_{2,1,1}=\Delta_{a\sigma_1}\text{ and hence}$$
$$g_{a,1}= x^{n-a\sigma_1}  \sum\limits_{p+q=a-1}T_{1,1}^pT_{1,\sigma_1+1}^q \pi_{1,p\sigma_1}
+ T_{1,1}^{a}\Delta_{a\sigma_1}.\tag\tnum{defg}$$Apply (\tref{defg}) with $a=a_1$ and $a=a_1+1$ and use (\tref{fact1}),  (\tref{fact2}), and  (\tref{fact3}), as in the proof of part (a), in order to establish (e). Likewise, set $a=f(\emptyset)$ in (\tref{defg}) and apply (\tref{fact4}) and (\tref{fact5}), as in the proof of part (d), to obtain (f). 
 
We prove (g) by showing that 
$$gT^{\pmb a}T_{k+1,1}^{f(\pmb a)}T_{k+1,j}\equiv y^n G_{(\pmb a,j)},$$whenever $\pmb a$ is an eligible $k$-tuple and $1\le j\le r(\pmb a)$.  Start with $\rho=2$. Recall that $g=h_{0,0}$ and $T_{3,1}=y$. If $\pmb a=(a_1,a_2)$ is an eligible tuple
and $1\le j\le r(\pmb a)$, then $r(\pmb a)=1=j$ and (a) and (b) show that 
$$gT_{1,1}^{a_1}T_{2,1}^{a_2}T_{3,1}^{f(\pmb a)}T_{3,j}\equiv y^{a_1\sigma_1+a_2\sigma_2+f(\pmb a)+1}h_{a_1,a_2}=y^nG_{(\pmb a,1)}.$$ If $(a_1)$ is an eligible tuple and $1\le j\le r(a_1)$, then (a), (b), and (c) yield 
$$g T_{1,1}^{a_1}T_{2,1}^{f(a_1)}T_{2,j}\equiv y^{a_1\sigma_1+f(a_1)\sigma_2+\sigma_2+1-r(a_1)}g_{a_1,j}=y^nG_{((a_1),j)}.$$ If $1\le j\le r(\emptyset)$, then (a) and (d) yield
$$g T_{1,1}^{f(\emptyset)}T_{1,j}\equiv y^{f(\emptyset)\sigma_1+\sigma_1+1-r(\emptyset)}f_{j}=y^nG_{(\emptyset,j)}.$$Now take $\rho=1$. Recall that $g=g_{0,1}$ and $y=T_{2,1}$. 
If $(a_1)$ is an eligible tuple and $1\le j\le r(a_1)$, then $r(a_1)=1=j$ and (e) gives 
$$gT_{1,1}^{a_1}T_{2,1}^{f(a_1)}T_{2,j}\equiv y^{a_1\sigma_1+f(a_1)+1}g_{a_1,1}=y^nG_{((a_1),1)}.$$ Finally, if $1\le j\le r(\emptyset)$, then (e) and (f) give
$$gT_{1,1}^{f(\emptyset)}T_{1,j}\equiv y^{f(\emptyset)\sigma_1+\sigma_1+1-r(\emptyset)}f_j=y^nG_{(\emptyset,j)}. \qed$$\enddemo

\remark{Remark \tnum{3.11}} Here is an alternative proof of Theorem \tref{A=L}. Let $t_1,\dots,t_{\rho}$ be variables and set $t_{\ell}=t_{\rho+1}=1$. Consider the homomorphism of $k[x,y]$-algebras
$$\alignat3 \pmb \pi\:k[\{T_{i,j}\mid 1\le i\le \ell, 1\le j\le \sigma_i+1\}]\phantom{,T_m}&\to k[x,y,t_1,t_{2}],&\quad&\text{if $\rho=2$, and}\hfill\\\vspace{5pt}
\pmb \pi\:k[\{T_{i,j}\mid 1\le i\le \ell, 1\le j\le \sigma_i+1\},T_m]&\to k[x,y,t_1,T_m],&\quad&\text{if $\rho=1$,}\hfill\endalignat$$
with $\pmb \pi(T_{i,j})=x^{j-1}y^{\sigma_i-j+1}t_i$, and, if $\rho=1$, with $\pmb \pi(T_m)=T_m$.  Notice that $\ker \pmb \pi=I_2(\psi)$; indeed, $I_2(\psi)$ is clearly contained in $\ker \pmb \pi$ and both sides are prime ideals of the same height $m-2$. Thus, to prove Theorem \tref{A=L}, it suffices to verify that
$$\pmb \pi(gT^{\pmb a}T_{k+1,1}^{f(\pmb a)}T_{k+1,j})=y^n \pmb \pi(G_{(\pmb a,j)})$$ for every $k$, $j$, and $\pmb a$, where $0\le k\le \rho$, $\pmb a$ is an eligible $k$-tuple, and $1\le j\le r(\pmb a)$. After evaluating the left hand side and dividing by $y^n$, the asserted equality becomes
$$\pmb \pi(G_{(\pmb a,j)})=\sum_{s=0}^nx^{n-s+j-1}y^{s+r(\pmb a)-j}t_1^{a_1}\cdots t_k^{a_k}t_{k+1}^{f(\pmb a)+1}\pmb \pi(c_s).$$From Definition \tref{1D5.3} one easily sees that $\pmb \pi(G_{(\pmb a,j)})$ is a sum of $n+1$ distinct terms of the form $m_s\pmb \pi(c_s)$, where $0\le s\le n$ and each $m_s$ is a monomial. By computing degrees in the various indeterminates, for instance, one deduces that indeed $m_s$ is equal to $x^{n-s+j-1}y^{s+r(\pmb a)-j}t_1^{a_1}\cdots t_k^{a_k}t_{k+1}^{f(\pmb a)+1}$. 
\endremark
\bigskip
\bigskip
\bigpagebreak
\SectionNumber=4\tNumber=1
\flushpar{\bf \number\SectionNumber.\quad Depth, reduction number, regularity, and Hilbert function.}
\medskip

This section is mainly about the special fiber ring $\Cal F(I)=\Cal R(I)/(x,y)$. We compute the depth, reduction number, and regularity of $\Cal F(I)$. A related invariant, the postulation number of $\Cal F(I)$, is computed in Corollary \tref{last}. Most of the results are collected in Theorem \tref{T1}; these results are proved, in a more general setting, in \cite{\rref{KPU-d}}; see Theorem \tref{4.3}. 
The main result of this section is Theorem \tref{red1} where we calculate the reduction number, $r(I)$, of $I$ when $\rho=2$. Observation \tref{ch} shows how we will use the rational normal scrolls of Section 1 to calculate $r(I)$. Theorem \tref{soc} is a general result connecting reduction number and Hilbert function for rings of minimal multiplicity; it is based on the Socle Lemma of Huneke-Ulrich. Proposition \tref{P39.2} is a curious result which allows us to circumvent the characteristic zero hypothesis in the Socle Lemma; we create a ring in which the bracket powers of the maximal ideal are equal to the ordinary powers, independent of the characteristic of the field. 

Let $B = \oplus_{i\ge 0}B_i$ be a standard graded Noetherian algebra
over an infinite field $k$ with $D$ equal to the Krull dimension of
$B$. The unique maximal homogenous ideal of $B$ is denoted by
$\maxm_B$. Let $I$ be an ideal of height $D$ generated by
homogeneous elements in $B$ of the same degree $d$. By a {\it
homogeneous minimal reduction} of $I$ we mean an ideal $J$ generated
by $D$ homogeneous elements in $I$ of degree $d$ so that
$I^{i+1}=JI^i$ for all large $i$.  Homogeneous minimal reductions exist; in fact any ideal generated by $D$ general $k$-linear combinations of forms of degree $d$ generating $I$ will do. The reduction number of $I$ with
respect to $J$ is
$$r_J(I)=\min\{i\ge 0\mid I^{i+1}=JI^i\},$$ and the {\it reduction
number} of $I$ is defined by $$r(I) =\min\{r_J(I)\mid \text{$J$ is a
homogenous minimal reduction of $I$}\}.$$

The homogenous minimal reductions of $\maxm_B$ are exactly the
ideals $J$ generated by linear systems of parameters, and the
reduction numbers can be characterized as
$$r_J(\maxm_B)=\min\{i\ge 0\mid \maxm_B^{i+1}\subset J\}.$$ Sometimes it is
convenient to write $r(B)$ in place of $r(\maxm_B)$. The reduction
number of an ideal $I$ is equal to be the reduction number of the
maximal homogeneous ideal in the special fiber ring $\Cal F(I)$;
that is $r(I)=r(\Cal F(I))$.

Recall that $B$ is said to have {\it minimal multiplicity} whenever $$e(B)=\operatorname{edim}(B)-\dim B+1.$$ This condition obtains if $r(B)\le 1$, in particular, if $\reg (B)\le 1$, and all three conditions are equivalent for a Cohen-Macaulay ring $B$. We will often use the fact that a one-dimensional Cohen-Macaulay ring has minimal multiplicity $e$ if and only if $H_{B}(i)=e$ for every $i\ge 1$. Standard examples of Cohen-Macaulay rings having minimal multiplicity include the algebras $A$ and $A\check{\phantom{a}}$ considered in Definition \tref{CN} and Data \tref{data4} below. 

The following notation is used often in this section.
\definition{Data \tnum{data4}}   Adopt the notation of Definition \tref{CN}.  Let $\check{\phantom{a}}$ mean image in $A\check{\phantom{a}}=A/(x,y)A$. 
\enddefinition

\proclaim{Observation \tnum{ch}} Adopt Data \tref{data4}. The following statements hold.
\smallskip \flushpar{\rm (1)} The ring $A\check{\phantom{a}}$ is defined by the maximal minors of  a scroll matrix and the ideal $K^{(n)}\check{\phantom{a}}$ is the $n^{\text{th}}$ symbolic power of a height one prime ideal of $A\check{\phantom{a}}$.
\smallskip \flushpar{\rm (2)} The special fiber ring of $I$ is equal to $\Cal F(I)=A\check{\phantom{a}}/\Cal A\check{\phantom{a}}$.
\smallskip \flushpar{\rm (3)} The graded $A\check{\phantom{a}}$-modules  $\Cal A\check{\phantom{a}}$ and $K^{(n)}\check{\phantom{a}}(-1)$ are isomorphic.
\smallskip \flushpar{\rm (4)} The reduction number of $I$ is equal to
$r(I)=r\({A\check{\phantom{a}}}/{K^{(n)}\check{\phantom{a}}}\)+1$.
 \endproclaim

\demo{Proof} Item (1)  is essentially obvious. The ring $A\check{\phantom{a}}$ equals $k[T_1,\dots,T_m]/I_2(\psi_{\text{\rm tr}})$, where  $\psi_{\text{\rm tr}}$ is the following truncation of $\psi$:
$$\psi_{\text{\rm tr}}=\cases \psi_1,&\text{if $\rho=1$},\\
\bmatrix \psi_1&\psi_2\endbmatrix,&\text{if $\rho=2$}.\endcases$$
A generating set of the ideal $K^{(n)}\check{\phantom{a}}$ of $A\check{\phantom{a}}$ is given in Theorem \tref{Kupn}.
On the other hand, one may consider the height one prime ideal $\kappa$ of $A\check{\phantom{a}}$ which is generated by the top row of $\psi_{\text{\rm tr}}$. A   generating set for the $n^{\text{th}}$ symbolic power, $\kappa^{(n)}$, of  $\kappa$ may also be found in Theorem \tref{Kupn}. The ideals 
  $K^{(n)}\check{\phantom{a}}$  and $\kappa^{(n)}$ of $A\check{\phantom{a}}$ have the same generators and therefore they are equal.  For (2), we have $\Cal F(I)=\Cal R(I)/(x,y)\Cal R(I)$ and $\Cal R(I)=A/\Cal AA$.

We prove (3). Recall from  (\tref{g}) that $$g=g(x,y)=\sum_{u=0}^nc_ux^{n-u}y^u\tag\tnum{gg}$$ and    from Theorem \tref{main1}(a) that the ideals
$$y^n\Cal AA \and g(x,y)K^{(n)}$$ of $A$ are equal. 
Fix a pair of subscripts $i,j$ with $1\le i\le \rho$ and $1\le j\le \sigma_i$.
Multiply both sides of (\tref{gg}) by $T_{i,j}^n$. Notice that, in $S$, $$\split T_{i,j}^ng(x,y)&{}=\sum_{u=0}^nc_u(T_{i,j}x)^{n-u}(T_{i,j}y)^u\equiv \sum_{u=0}^nc_u(T_{i,j+1}y)^{n-u}(T_{i,j}y)^u\\&{}=y^ng(T_{i,j+1},T_{i,j}).\endsplit $$ (See Remark \tref{1R6.8}(d) for the meaning of $\equiv$.)  Conclude that 
$$y^nT_{i,j}^n\Cal AA =y^n g(T_{i,j+1},T_{i,j})K^{(n)}.$$ The ring $A$ is a domain and $y\neq 0$; so, $$T_{i,j}^n\Cal AA=g(T_{i,j+1},T_{i,j})K^{(n)}. $$ 
In particular, we have $$T_{1,1}^n\Cal A\check{\phantom{a}}=g(T_{1,2},T_{1,1})K^{(n)}\check{\phantom{a}}.$$Now $T_{1,1}$ has non-zero image in the domain $A\check{\phantom{a}}$, and so does 
$\Cal A$ because $$\dim \Cal F(I)=2<3=\dim   A\check{\phantom{a}}.$$ Thus, the image of 
$g(T_{1,2},T_{1,1})$ in $A\check{\phantom{a}}$ cannot be zero either. It follows that both 
$T_{1,1}^n$ and $g(T_{1,2},T_{1,1})$ are non zerodivisors on the domain $A\check{\phantom{a}}$.
 Assertion (3) is proved. 

We prove (4). We have seen that 
$$r(I)=r(\Cal F(I))=r\(  {A\check{\phantom{x}}}/{\Cal A\check{\phantom{x}}}\).$$
The reduction numbers of the two-dimensional standard graded rings   $A\check{\phantom{x}}/{\Cal A\check{\phantom{x}}}$ and 
$A\check{\phantom{x}}/K^{(n)}\check{\phantom{x}}$ 
may be computed after reducing modulo two generic linear forms,
in which case the reduction number is simply the top socle degree, 
 see,  \cite{\rref{SUV}, Lemma~3.4}. Let $k(u)$ be the appropriate purely transcendental extension of $k$, let $\ell_1$ and $\ell_2$ be two generic linear forms in $A\check{\phantom{x}}\t_kk(u)$, and let $\overline{\phantom{x}}$ represent image in 
$\overline{A\check{\phantom{x}}}=(A\check{\phantom{x}}\t_kk(u))/(\ell_1,\ell_2)(A\check{\phantom{x}}\t_kk(u))$. For every non-zero homogeneous element $z$ in $A\check{\phantom{x}}$ of positive degree, the sequence
 $z,\ell_1,\ell_2$ is a regular  on ${A\check{\phantom{x}}\t_kk(u)}$. Thus, both $T_{1,1}^n$ and $g(T_{1,2},T_{1,1})$ are non zerodivisors   on ${A\check{\phantom{x}}\t_kk(u)}$. It follows that 
the graded $\overline{A\check{\phantom{x}}}$-modules  $\overline{\Cal A\check{\phantom{x}}}$ and $\overline{K^{(n)}\check{\phantom{x}}}(-1)$    are isomorphic.
We know, since $\overline{A\check{\phantom{x}}}$ has minimal multiplicity, that $\lambda\(\overline{A\check{\phantom{x}}}_s\)=m-2$ for all $s\ge  1$ (see also  Proposition \tref{P11.4}); so,
$$\align \(\overline{A\check{\phantom{x}}}/\overline{\Cal A\check{\phantom{x}}}\)_s=0
&{}\iff \lambda\(\overline{\Cal A\check{\phantom{x}}}\)_s=m-2
\iff \lambda\(\overline{K^{(n)}\check{\phantom{x}}}\)_{s-1}=m-2\\&{}
\iff  \(\overline{A\check{\phantom{x}}}/\overline{K^{(n)}\check{\phantom{x}}}\)_{s-1}=0, \quad\text{and}\endalign $$ 
$$\align r(I)=r\(A\check{\phantom{x}}/\Cal A\check{\phantom{x}}\)&{}=\max\{s\mid \(\overline{A\check{\phantom{x}}}/\overline{\Cal A\check{\phantom{x}}}\)_s\neq 0\}\\&{}=\max\{s\mid \(\overline{A\check{\phantom{x}}}/\overline{K^{(n)}\check{\phantom{x}}}\)_s\neq 0\}+1=r\(A\check{\phantom{x}}/K^{(n)}\check{\phantom{x}}\)+1. \qed \endalign$$
\enddemo
 
\proclaim{Theorem \tnum{T1}} Adopt Data \tref{data4} with $n\ge 2$. 
\medskip\flushpar{\rm (a)} We have $$\depth {\operatorname{gr}}_I(R) + 1=\depth \Cal R(I)=\depth \Cal F(I)=\cases 2,&\text{if $\rho=1$},\\1,&\text{if $\rho=2$.}\endcases$$ In particular, $\Cal F(I)$ is Cohen-Macaulay if and only if $\rho=1$. 
\medskip\flushpar{\rm (b)} If $\rho=1$, then  $$\tsize r(I)=r(\Cal F(I))=\reg
\Cal F(I)=\lceil\frac{n-1}{\sigma_1} \rceil + 1.$$
\medskip\flushpar{\rm (c)} If $\rho=2$, then $\reg (\Cal F(I))=\lceil\frac{n-1}{\sigma_2}\rceil+1$. 
\endproclaim\remark{Remark} The value of $r(I)=r(\Cal F(I))$ when $\rho=2$ is computed in Theorem \tref{red1}.\endremark
\demo{Proof}
We apply Theorem \tref{main1} and Observation \tref{ch}. The isomorphisms $$\tsize\Cal R(I)=\frac{A}{\Cal AA},\quad \Cal F(I)=\frac{A\check{\phantom{x}}}{\Cal A\check{\phantom{x}}},\quad \Cal AA\iso K^{(n)}(-1),\and \Cal A\check{\phantom{x}}\iso K^{(n)}\check{\phantom{x}}(-1)$$ tell us that 
$$\depth \Cal R(I)=\depth A/K^{(n)}\and \depth \Cal F(I)= \depth A\check{\phantom{x}}/K^{(n)}
\check{\phantom{x}}.$$
They also show that $\reg \Cal A\check{\phantom{x}}=\reg K^{(n)}\check{\phantom{x}}+1$. As $0\neq K^{(n)}\check{\phantom{x}} \subsetneq A\check{\phantom{x}}$, we have $\reg K^{(n)}\check{\phantom{x}}\ge 1$, and therefore $\reg \Cal A\check{\phantom{x}}\ge 2$. On the other hand, $\reg  A\check{\phantom{x}}=1$; hence, $\reg \Cal A\check{\phantom{x}}> \reg   A\check{\phantom{x}}$. This strict inequality gives $\reg (A\check{\phantom{x}}/\Cal A\check{\phantom{x}})=\reg \Cal A\check{\phantom{x}}-1$. It follows that 
$$ \reg \Cal F(I)= \reg  (A\check{\phantom{x}}/\Cal A\check{\phantom{x}})= \reg (\Cal A\check{\phantom{x}})-1=\reg K^{(n)}\check{\phantom{x}}.$$
One may now read $\reg \Cal F(I)=\lceil\frac{n-1}{\sigma_{\rho}}\rceil+1$ from Theorem \tref{4.3}. If $\rho=2$, then one may also read $\depth \Cal R(I)=\depth \Cal F(I)=1$. 
If $\rho=1$, then the variable  $T_m$ is 
not involved in
$$\psi=\bmatrix T_1&\dots&T_{m-2}&y\\T_2&\dots&T_{m-1}&x\endbmatrix, \quad
\psi_{\text{\rm tr}}=\bmatrix T_1&\dots&T_{m-2}\\T_2&\dots&T_{m-1}\endbmatrix,$$
$K^{(n)}$, or $K^{(n)}\check{\phantom{x}}$; so $T_m$ is regular on $A/K^{(n)}$ and 
 $A\check{\phantom{x}}/K^{(n)}\check{\phantom{x}}$ and Theorem \tref{4.3} tells us that $$\depth A/(T_m,K^{(n)})
=\depth A\check{\phantom{x}}/(T_m,K^{(n)}\check{\phantom{x}})=1.$$ Therefore $\depth \Cal R(I)=\depth \Cal F(I)=2$ in this case. For any value of $\rho$, the Rees algebra $\Cal R(I)$ is not Cohen-Macaulay. It follows that ${\operatorname{gr}}_I(R)$ is not Cohen-Macaulay either by \cite{\rref{L}, Thm\.~5} and then $\depth \Cal R(I) =  \depth {\operatorname{gr}}_I(R) + 1$ by \cite{\rref{HM}, Thm\.~3.10}. Finally, we recall that 
 if $\Cal F(I)$ is a Cohen-Macaulay ring, then $\reg \Cal F(I)=r(I)$; indeed both quantities are equal to the top socle degree of $\Cal F(I)$ modulo a  linear system of parameters. \qed \enddemo

\proclaim{Theorem \tnum{4.3}}Let $\sigma_1\ge \dots\ge \sigma_{\ell}\ge 1$ and   $n\ge 2$ be integers, and let
 $P$ be the polynomial ring $$k[\{T_{i,j}\mid 1\le i\le \ell\and 1\le j \le \sigma_i+1\}].$$
For each $u$, with $1\le u\le \ell$, let $\psi_u$ be 
the generic scroll matrix $$\psi_u=\bmatrix T_{u,1}&T_{u,2}&\dots&T_{u,\sigma_u-1}&T_{u,\sigma_u}\\
T_{u,2}&T_{u,3}&\dots&T_{u,\sigma_u}&T_{u,\sigma_u+1}\endbmatrix. $$  Let $\Psi$  be the matrix
$$\Psi=\bmatrix \psi_1&\vrule&\dots&\vrule&\psi_{\ell}\endbmatrix, $$ 
    $A$ be the ring $P/I_2(\Psi)$, and $K$ be the ideal in $A$ generated by the entries of the top row of $\Psi$. 
Then
$$\tsize \depth A/K^{(n)}=1\and \reg K^{(n)}=\lceil \frac{n-1}{\sigma_{\ell}}
 \rceil + 1.$$
\endproclaim
\demo{Proof} See \cite{\rref{KPU-d}, Theorem 5.5 and Corollary 2.6}. \qed \enddemo

\proclaim{Theorem \tnum{red1}}Adopt Data \tref{data4}. If $I$ is as in Definition \tref{CN} with $\rho=2$, then the following statements hold.
\medskip\flushpar{\rm(a)} $\tsize
\lceil\frac{n}{\sigma_1} \rceil \leq r(I) \leq
\lceil\frac{n-1}{\sigma_1} \rceil + 1$.
\medskip\flushpar{\rm(b)} 
If $\sigma_1 | n-1$, then
$r(I)=\lceil\frac{n}{\sigma_1} \rceil =
\lceil\frac{n-1}{\sigma_1} \rceil + 1$.
\medskip\flushpar{\rm(c)} $\tsize r(I)=\lceil\frac{n}{\sigma_1} \rceil\iff \Hgy_
{K^{(n)}\check{\phantom{x}}}(\lceil\frac{n}{\sigma_1} \rceil) \geq m-2$.
\endproclaim
\remark{Remarks} 
1. The exact value of $\Hgy_
{K^{(n)}\check{\phantom{x}}}(\lceil\frac{n}{\sigma_1} \rceil)$ depends on the interaction between the three integers $\sigma_1$, $\sigma_2$, and $n$, and is not 
difficult to calculate.  From Theorem \tref{Kupn} we know that $K^{(n)}\check{\phantom{x}}$ is equal to
$$\eightpoint \cases T_{1,1}^{f(\emptyset)}(T_{1,1},\dots, T_{1,r(\emptyset)})A\check{\phantom{a}},&\text{if $\rho=1$},\\
T_{1,1}^{f(\emptyset)}(T_{1,1},\dots, T_{1,r(\emptyset)})A\check{\phantom{a}}
+\sum\limits_{a_1=0}^{f(\emptyset)}T_{1,1}^{a_1}T_{2,1}^{f(a_{1})}(T_{2,1},\dots, T_{2,r(a_1)})A\check{\phantom{a}},
&\text{if $\rho=2$}.\endcases\tag\tnum{4.55}$$
Furthermore, the generators listed here form a homogeneous minimal generating set of 
$K^{(n)}\check{\phantom{x}}$  according to \cite{\rref{KPU-d}, Prop.~1.20}. Clearly, $K^{(n)}\check{\phantom{x}}$ is concentrated in degrees at least $f(\emptyset)+1=\lceil \frac{n}{\sigma_1}\rceil$, Thus, to compute $\Hgy_
{K^{(n)}\check{\phantom{x}}}(\lceil\frac{n}{\sigma_1} \rceil) $ one only needs to count the minimal homogeneous generators of degree $\lceil \frac{n}{\sigma_1}\rceil$. Thus, the exact value of  $\Hgy_
{K^{(n)}\check{\phantom{x}}}(\lceil\frac{n}{\sigma_1} \rceil)$ is equal to 
$\sigma_1\lceil\frac{n}{\sigma_1} \rceil-n+1$ plus
the cardinality of the following set
$$\tsize \{(i,j,k)\mid i+j+1= \lceil\frac{n}{\sigma_1} \rceil \quad\text{and}\quad \sigma_1i+\sigma_2j<n\le \sigma_1i+\sigma_2(j+1)+1-k\},$$where $i$ and $j$  are non-negative integers and $k$ is a positive integer.

\medskip\flushpar 2. We prove (a) now. Assertion (b) is obvious. 

\medskip\flushpar 3. Part (a)  shows that there are only two possible choices for $r(I)$. Furthermore, in the proof of (a),  we   learn a necessary condition for $r(I)$ to take on the smaller of the two values; namely, that $K^{(n)}\check{\phantom{x}}$ contain at least $m-2$ linearly independent homogeneous elements of degree $\lceil\frac n{\sigma_1}\rceil$. The proof that this condition is sufficient (i.e., part (c)) appears at the end of this section.  \endremark

\demo{Proof of {\rm(a)}}Use the notation of Observation \tref{ch}. So, $A\check{\phantom{x}}=k[T_1,\dots,T_m]/I_2(\psi_{\text{\rm tr}})$ and $\Cal F(I)=A\check{\phantom{x}}/\Cal A\check{\phantom{x}}$. The convention of (\tref{Alt}) is in effect and each variable $T_i$ has two names. 

We first establish the inequality on the right.
Let $J= (T_m,
T_{\sigma_1+1}-T_{\sigma_1 + 2})$ and write $^{-}$ for images in
$\overline{A\check{\phantom{x}}}=A\check{\phantom{x}}/JA\check{\phantom{x}}$.
Notice that $T_m$,
$T_{\sigma_1+1}-T_{\sigma_1 + 2}$ form a linear system of parameters in
$A\check{\phantom{x}}/K^{(n)}\check{\phantom{x}}$. 
In conjunction with Observation \tref{ch}(4), this yields 
$$r(I) =r\({A\check{\phantom{a}}}/{K^{(n)}\check{\phantom{a}}}\)+1 
 \le \max \left\{i \left\vert \[\overline{A\check{\phantom{x}}}/\overline{K^{(n)}\check{\phantom{x}}}\]_i
\not=0\right.\right\}+1.$$ Clearly, $[\overline{A\check{\phantom{x}}}/(T_{1,1})]_i=0$, for $i>1$; hence, 
$\overline{A\check{\phantom{x}}}_i=T_{1,1}^{i-1}\overline{A\check{\phantom{x}}}_1$ for $i\ge 1$. On the other hand,
the generators of $\overline{K^{(n)}\check{\phantom{x}}}$ are listed in (\tref{4.55}). Observe that $T_{2,1}T_{2,j}=0$ in $\overline{A\check{\phantom{x}}}$, for $1\le j\le \sigma_2+1$, because
$$T_{2,1}T_{2,j}=T_{1,\sigma_1+1}T_{2,j}=T_{1,\sigma_1-\sigma_2+j}T_{2,\sigma_2+1}=0.$$
Notice that all variables make sense. The first and last equations are due to $J$. The middle equation  happens because of the determinantal relations. So, $$\overline{K^{(n)}\check{\phantom{x}}}=T_{1,1}^{f(\emptyset)}\[(T_{1,1},\dots,T_{1,r(\emptyset)})+\chi(f(f(\emptyset))=0)(T_{2,1},\dots,T_{2,r(f(\emptyset))})\].\tag\tnum{4.5}$$ Observe immediately that 
$$\overline{A\check{\phantom{x}}}_i=T_{1,1}^{i-1}\overline{A\check{\phantom{x}}}_1\subseteq \overline{K^{(n)}\check{\phantom{x}}}\tag\tnum{4.6}$$ for $i>f(\emptyset)+1$. Recall that
$$r(\emptyset)=\sigma_1\iff \sigma_1|(n-1)\iff f(\emptyset)\sigma_1=n-1.\tag\tnum{4.7}$$
If (\tref{4.7}) occurs, then $f(f(\emptyset))=0$ and $r(f(\emptyset))=\sigma_2$. Therefore, if 
(\tref{4.7}) occurs, then (\tref{4.5}) shows that (\tref{4.6}) also occurs at the value $i=
f(\emptyset)+1$. 
We  see that  
$$\tsize r(I)-1\le\max \left\{i \left\vert \[\overline{A\check{\phantom{x}}}/\overline{K^{(n)}\check{\phantom{x}}}\]_i\neq 0\right.\right\}
\le\cases
f(\emptyset),&\text{if $\sigma_1|(n-1)$,}\\f(\emptyset)+1,&\text{if $\sigma_1\not|(n-1)$;}\endcases$$ hence, $r(I)-1\le  \lceil\frac {n-1}{\sigma_1}\rceil$.
 
Now we establish the inequality on the left. We may assume that the field is infinite. Let $\ell_1$ and $\ell_2$ be two general linear forms in $k[T_1,\dots,T_m]$, $J$ be the ideal $(\ell_1,\ell_2)$ of $A\check{\phantom{x}}$, and $\overline{\phantom{x}}$ represent image in $\overline{A\check{\phantom{x}}}=A\check{\phantom{x}}/JA\check{\phantom{x}}$. We see that $\ell_1,\ell_2$ is a general linear system of parameters in $A\check{\phantom{x}}/K^{(n)}\check{\phantom{x}}$; hence, 
$$r\(A\check{\phantom{x}}/K^{(n)}\check{\phantom{x}}\)=\max \left\{i \left\vert\[\overline{A\check{\phantom{x}}}/\overline{K^{(n)}\check{\phantom{x}}}\]_i\neq 0\right.\right\},$$see, for instance \cite{\rref{T03}, Cor\.~2.2}.
However, $$\[\overline{A\check{\phantom{x}}}/\overline{K^{(n)}\check{\phantom{x}}}\]_{f(\emptyset)}\neq 0$$ because $\overline{A\check{\phantom{x}}}$ is a ring of positive Krull dimension and the minimal generator degree of $\overline{K^{(n)}\check{\phantom{x}}}$ is greater than $f(\emptyset)$, see (\tref{4.55}). Thus, $$ f(\emptyset)\le r\(A\check{\phantom{x}}/K^{(n)}\check{\phantom{x}}\)=r(I)-1 $$
and then $\lceil\frac{n}{\sigma_1} \rceil=f(\emptyset)+1\le r(I)$.
\qed 
\enddemo

The proof of Theorem {\rm\tref{red1}}{\rm(c)} will be based on the next general result
relating reduction numbers to Hilbert functions. Assume that $(B, \maxm)$ is a one-dimensional Cohen-Macaulay standard graded ring over a field. Assume also that $B$ has 
  minimal multiplicity $e$ and $L$ is a homogeneous $\maxm$-primary ideal. In this case  $ r(B/L)< s$ if and only if $\Hgy_L(s)\geq e \, .$
The purpose of the following theorem is to prove an analogous statement in dimension two. 
The first difference function of the Hilbert function is denoted $\Delta\Hgy$ and is defined by $\Delta\Hgy_M(i)=\Hgy_M(i)-\Hgy_M(i-1)$. 

\proclaim{Theorem \tnum{soc}}
Let $(B,\maxm)$ be a standard graded domain over a field $k$.
Assume that
$B$ is a two-dimensional Cohen-Macaulay ring with
minimal multiplicity $e$. Let $s$ be a positive integer and $L$ be a homogeneous ideal of $B$ with $\dim B/L=1$. Let $\ell$ be a generic linear form defined over a purely transcendental extension field $k'$ of $k$ 
and assume that, after making a possible further field extension $k''$ of $k'$, 
$$\maxm^s \subset (\{y^s\mid y\in B_1 \otimes_k k''\})+(L,\ell).\tag\tnum{cnd}$$  
One has
$$ r(B/L)< s \Longleftrightarrow \Delta \Hgy_L(s)\geq e.
$$
\endproclaim
\remark{Remark} The hypothesis (\tref{cnd}) is automatically satisfied if the characteristic of $k$ is zero.\endremark
\demo{Proof}  We replace $B$ by $B \otimes_k k'$. This does not change the reduction number of $B/L$ according to
\cite{\rref{SUV}, Lemma~3.4}. Write $^{-}$ for images in $\overline{B}=B/(\ell)$. Notice that
$\overline{B}$ is a domain by \cite{\rref{Ho}, Theorem} and
$r(B/L)=r\(\overline{B}/\overline{L}\)$ again by
\cite{\rref{SUV}, Lemma~3.4}, whereas  $\Delta \Hgy_L= \Hgy_{L/\ell L}$.  Moreover,  $\overline{B}/\overline{L}$ is an Artinian ring. We need to show that
$\overline{L}_{s}=\overline{B}_{s}$ if and only if $\Hgy_{L/\ell
L}(s)\geq e$. Since $\Hgy_{\overline{B}}(s)=e$ it remains to prove that
$\Hgy_{\overline{L}}(s)\geq e$ if and only if $\Hgy_{L/\ell L}(s) \geq e$. As
$\Hgy_{\overline{L}}(s) \leq \Hgy_{L/\ell L}(s)$ it suffices to show that
if $\Hgy_{L/\ell L}(s) \geq e$ then $\Hgy_{\overline{L}}(s) \geq e$.

If $\overline{B}/\overline{L}$ has a non-zero homogeneous socle
element of degree $j < s$, then $z \overline{\m} \subset
\overline{L}$ for some homogeneous non-zero element $z \in
\overline{B}$ of degree $j$. Hence $$\Hgy_{\overline{L}}(s) \geq
\Hgy_{z\overline{\m}}(s) = \Hgy_{\overline{\m}}(s-j)$$ because
$\overline{B}$ is a domain. Clearly $\Hgy_{\overline{\m}}(s-j) = e$
since $s-j \geq 1$, which gives $\Hgy_{\overline{L}}(s) \geq e$.

Thus we may assume that the socle of $\overline{B}/\overline{L}$  is concentrated in degrees $\geq s$. For the remainder of the proof we do not need anymore that $\overline{B}$ is a domain, thus we may extend the ground field to assume that condition (\tref{cnd}) holds. 
We wish to apply the Socle Lemma  \cite{\rref{HU}, Cor\.~3.11(i)} to the exact sequence 
$$
0 \longrightarrow (L \link \ell/L)(-1)\longrightarrow (B/L)(-1)
@> \ell >> B/L @>>>
B/(L,\ell)=\overline{B}/\overline{L} \longrightarrow 0 \, .
$$
The statement of the Socle Lemma requires that the field have characteristic zero; however, this hypothesis is only used in order to ensure that condition (\tref{cnd}) is satisfied. Also, the Socle
Lemma requires $\ell$ to be a general linear form, but the proof
also works for generic linear forms, see \cite{\rref{HU}, Prop\.~3.5}. So the Socle Lemma may be applied in the present situation. In our
setting it says that if the socle of $\overline{B}/\overline{L}$ is
concentrated in degrees $\geq s$ then so is $L \link \ell/L$. Thus 
$[L \link \ell]_{s-1}=L_{s-1}$. It follows that
$$[L \cap (\ell)]_s= [\ell(L \link \ell)]_s=\ell[L \link \ell]_{s-1}=\ell L_{s-1}=
[\ell
L]_s.$$
This gives $[\overline{L}]_s=[L/(\ell L)]_s$. Thus
$\Hgy_{\overline{L}}(s)= \Hgy_{L/\ell L}(s) \geq e$. \qed\enddemo

In the next proposition we show that the homogeneous coordinate ring of any rational normal curve satisfies assumption  (\tref{cnd}) regardless of the characteristic.
If one is only interested in characteristic zero, then Proposition \tref{P39.2} may be skipped. 
\proclaim{Proposition \tnum{P39.2}} Let $k'$ be the field of rational functions $k(\lambda_1,\dots,\lambda_c)$ in $c$ variables over the field $k$ and let $k''$ be any extension field of $k'$ for which the polynomial
$$P(x)=x^c-\lambda_cx^{c-1}-\dots-\lambda_2 x-\lambda_1$$ splits into linear factors.  Let $A$ be the standard graded algebra $k''[T_1,\dots,T_c]/I_2(\psi)$, where $\psi$ is the scroll matrix 
$$\psi=\bmatrix T_1&T_2&\dots&T_{c-1}&T_c\\T_2&T_3&\dots&T_{c}&\sum\limits_{i=1}^c\lambda_{i}T_{i}\endbmatrix.$$Then there exist homogeneous linear forms $v_1,\dots,v_c$ in $A_1$ such that 
$v_1^s,\dots,v_c^s$ is a $k''$-basis for $A_s$ for all $s\ge 1$.\endproclaim

\demo{Proof} Recall that $A$ is a one-dimensional standard graded ring of minimal multiplicity $c$ and $T_1$ is an $A$-regular linear form. (Some readers will find Proposition \tref{P11.4} to be helpful at this point.) Hence for any $s\geq 1$ one has  
$\dim A_s=c$ and $T_1^{s-1}T_1,\dots,T_1^{s-1}T_c$ form a basis of $A_s$. 
Once we have identified suitable candidates for $v_1,\dots,v_c$, then we need only verify that  $v_1^s,\dots,v_c^s$ are linearly independent. 
Ultimately, we pick $v_1,\dots,v_c$ to be a basis for $A_1$ which yields a simultaneous diagonalization of all of the linear transformations $\varphi_j=T_1^{-1}T_j\:A_1\to A_1$. 

Let $k'''\subseteq k''$ be 
the splitting field of $P(x)$ over $k'$. Since $\lambda_1, \ldots, \lambda_c$ are variables over $k$, the polynomial $P(x)$ has $c$ distinct roots. In particular, the field extension $k' \subset k'''$ is separable. 

The matrix representation of the endomorphism $\varphi_2= T_1^{-1}T_2\: A_1 \to A_1$ with respect to the basis $T_1, \ldots, T_c$ is 
 
$$\bmatrix
0& & & & & & \lambda_1\\
1&0 & \cdot  & & & & \cdot\\
&1 & \cdot & \cdot & & & \cdot \\
& & &  \cdot & \cdot & & \cdot\\
& & & & \cdot & 0& \cdot \\
& & & & &1 &\lambda_c 
\endbmatrix.$$
This being a companion matrix it follows that the minimal polynomial of $\varphi_2$ is $P(x)$, which has $c$ distinct roots in $k'''$. Thus $\varphi_2$ is diagonalizable over $k'''$ with eigenvectors,  say, $v_1, \ldots, v_c$. 

On the other hand, for $2\le j \le c$ one has  $T_1T_j=T_2T_{j-1}$, hence $\varphi_j=\varphi_2\varphi_{j-1}$. Thus one sees by induction on $j$ that $v_1, \ldots, v_c$ are eigenvectors 
for every $\varphi_j$. In fact, there exist $\alpha_{i} \in k'''$ with 
$$T_j v_i= \alpha_{i}^{j-1}T_1v_i  \, .$$
Thus, $\maxm v_i \subset AT_1v_i$ for $1 \le i \le c$, and then  ${\maxm}^{s-1} v_i \subset AT_1^{s-1}v_i \,.$
In particular, $v_i^{s} \in AT_1^{s-1}v_i \, ,$ say 
$$
v_i^s=\beta_iT_1^{s-1}v_i, \tag\tnum{eig}
$$
for some $\beta_i \in k'''$.

Recall that $k'[T_1,\ldots, T_c]/I_2(\Psi)$ is a domain and $k' \subset k'''$ is a separable field extension. Therefore $k'''[T_1,\ldots, T_c]/I_2(\Psi)$ is reduced, hence $v_i^s \not=0$
in $$k'''[T_1,\ldots, T_c]/I_2(\Psi)\subseteq A,$$
and  $\beta_i \not=0$. Since 
$T_1^{s-1}$ is a non zerodivisor and $\beta_i$ are non-zero scalars,   (\tref{eig})  shows that 
$v_1^s, \ldots, v_c^s$ are indeed linearly independent over $k'''$ and therefore also over $k''$.\qed\enddemo

\demo{Proof of part {\rm(c)} of Theorem {\rm\tref{red1}}} Start with the ring $A\check{\phantom{x}}$ and the ideal $K^{(n)}\check{\phantom{x}}$ of Observation \tref{ch}. 
From Theorem \tref{red1}(a) and Observation \tref{ch}(4) we know that $$\tsize \lceil\frac n{\sigma_1}\rceil\le r(I)=r\(A\check{\phantom{x}}/K^{(n)}\check{\phantom{x}}\)+1 \, .$$
Hence it suffices to show that $$\tsize r\(A\check{\phantom{x}}/K^{(n)}\check{\phantom{x}}\)<\lceil\frac{n}{\sigma_1} \rceil\iff 
\Hgy_ {K^{(n)}\check{\phantom{x}}}(\lceil\frac{n}{\sigma_1} \rceil)\geq m-2.$$ 

First assume that $\Hgy_ {K^{(n)}\check{\phantom{x}}}(\lceil\frac{n}{\sigma_1} \rceil)\geq m-2$. Let $\ell$ be the linear form $T_{1, \sigma_1+1}-T_{2,1}$ of $A\check{\phantom{x}}$,  and let $\overline{\phantom{x}}$ denote images in the ring $\overline{A\check{\phantom{x}}}=A\check{\phantom{x}}/(\ell)\, .$  Clearly, $\ell$ is a regular element on $A\check{\phantom{x}}$. 
Notice that the image of $T_{2,1}$ in $A\check{\phantom{x}}$ lies in $K\check{\phantom{x}}$, whereas the image of $T_{1,\sigma_1+1}$ does not.
Therefore $\ell= T_{1, \sigma_1+1}-T_{2,1}$ is regular on $A\check{\phantom{x}}/(K\check{\phantom{x}})^{(n)}=A\check{\phantom{x}}/K^{(n)}\check{\phantom{x}}$. It follows that $\Hgy_ {\overline{K^{(n)}\check{\phantom{x}}}}(\lceil\frac{n}{\sigma_1} \rceil) = \Delta \Hgy_ {K^{(n)}\check{\phantom{x}}}(\lceil\frac{n}{\sigma_1} \rceil) \, .$ However, $K^{(n)}\check{\phantom{x}}$ is concentrated in degrees at least $\lceil\frac{n}{\sigma_1} \rceil$ and therefore $ \Delta \Hgy_ {K^{(n)}\check{\phantom{x}}}(\lceil\frac{n}{\sigma_1} \rceil)=   \Hgy_ {K^{(n)}\check{\phantom{x}}}(\lceil\frac{n}{\sigma_1} \rceil) \, . $ On the other hand,
$r\(\overline{A\check{\phantom{x}}}/\overline{K^{(n)}\check{\phantom{x}}}\) \ge r\(A\check{\phantom{x}}/K^{(n)}\check{\phantom{x}}\)$. Hence, it suffices to prove that  
$$\tsize \Hgy_ {\overline{K^{(n)}\check{\phantom{x}}}}(\lceil\frac{n}{\sigma_1} \rceil) \geq m-2\implies  r\(\overline{A\check{\phantom{x}}}/\overline{K^{(n)}\check{\phantom{x}}}\) <\lceil\frac{n}{\sigma_1} \rceil.$$
For this we wish to apply Theorem \tref{soc} to the integer $\lceil\frac{n}{\sigma_1} \rceil$ and the ideal  $\overline{K^{(n)}\check{\phantom{x}}}$ of the ring 
$\overline{A\check{\phantom{x}}}$. Notice that $\overline{A\check{\phantom{x}}}$ is the homogeneous coordinate ring of a rational normal curve. In particular, it is a two dimensional Cohen-Macaulay domain with minimal multiplicity $m-2$. By Proposition \tref{P39.2} the ring  $\overline{A\check{\phantom{x}}}$ satisfies condition (\tref{cnd}). Furthermore, 
$\overline{K^{(n)}\check{\phantom{x}}}$ is a  homogeneous ideal with
$\dim \overline{A\check{\phantom{x}}}/\overline{K^{(n)}\check{\phantom{x}}}=1$; thus, Theorem \tref{soc} implies that  $r\(\overline{A\check{\phantom{x}}}/\overline{K^{(n)}\check{\phantom{x}}}\) <\lceil\frac{n}{\sigma_1} \rceil$ if $\Delta \Hgy_ {\overline{K^{(n)}\check{\phantom{x}}}}(\lceil\frac{n}{\sigma_1} \rceil) \geq m-2\, .$ But again, 
$\Delta \Hgy_ {\overline{K^{(n)}\check{\phantom{x}}}}(\lceil\frac{n}{\sigma_1} \rceil) =  \Hgy_ {\overline{K^{(n)}\check{\phantom{x}}}}(\lceil\frac{n}{\sigma_1} \rceil) $. This completes the proof of the first implication.

Conversely, assume that $r\(A\check{\phantom{x}}/K^{(n)}\check{\phantom{x}}\)<\lceil\frac{n}{\sigma_1} \rceil$. Now let  $\overline{A\check{\phantom{x}}}$ denote the ring obtained from 
$A\check{\phantom{x}}$ by a purely transcendental extension of the field $k$ and by factoring out two generic linear forms. Write   $\overline{K^{(n)}\check{\phantom{x}}} = 
K^{(n)}\check{\phantom{x}} \overline{A\check{\phantom{x}}}$.  One has 
$$r\(\overline{A\check{\phantom{x}}}/\overline{K^{(n)}\check{\phantom{x}}}\) = r\(A\check{\phantom{x}}/K^{(n)}\check{\phantom{x}} \),$$ see \cite{\rref{SUV}, Lemma~3.4}. Therefore $r\(\overline{A\check{\phantom{x}}}/\overline{K^{(n)}\check{\phantom{x}}}\)<\lceil\frac{n}{\sigma_1} \rceil$. Because $\overline{A\check{\phantom{x}}}/\overline{K^{(n)}\check{\phantom{x}}}$
is Artinian and $\overline{A\check{\phantom{x}}}$ is a one-dimensional standard graded Cohen-Macaulay ring with minimal multiplicity $m-2$, we conclude that 
$ \Hgy_ {\overline{K^{(n)}\check{\phantom{x}}}}(\lceil\frac{n}{\sigma_1} \rceil) =   \Hgy_{\overline{A\check{\phantom{x}}}} (\lceil\frac{n}{\sigma_1} \rceil)=m-2\, .$ Clearly, 
$$\tsize \Hgy_ {K^{(n)}\check{\phantom{x}}}(\lceil\frac{n}{\sigma_1} \rceil) \geq \Hgy_ {\overline{K^{(n)}\check{\phantom{x}}}}(\lceil\frac{n}{\sigma_1} \rceil).$$ Hence indeed $\Hgy_ {K^{(n)}\check{\phantom{x}}}(\lceil\frac{n}{\sigma_1} \rceil) \geq m-2$.  \qed\enddemo

\bigskip
\bigskip
\bigpagebreak
\SectionNumber=5\tNumber=1
\flushpar{\bf \number\SectionNumber.\quad  Generalized Eagon-Northcott modules.}
\medskip
Let $I$ be the ideal of (\tref{Goal}). In Theorem \tref{pwrs4} we record the graded Betti numbers of $I^s$ for all $s$. The main step in the proof of this theorem is the calculation of the Hilbert function of $I^s$ and we do this by calculating 
$\lambda((S/H)_{(u,s)})$ and $\lambda(K^{(n)}_{(u,s)})$ for each bi-degree $(u,s)$. The $S$-module $S/H$ is resolved by an Eagon-Northcott complex and we have identified a filtration $\{\Cal E_{\pmb a}\}$ of $K^{(n)}$ so that each factor $\Cal E_{\pmb a}/\Cal D_{\pmb a}$ is a ``Generalized Eagon-Northcott module'', in the sense that it is resolved  by a generalized Eagon-Northcott complex.  See \cite{\rref{BV88}, Section 2C} or \cite{\rref{E95}, Section A2.6} for more information about these modules and complexes.
We define the generalized Eagon-Northcott modules in Definition \tref{D10.5}. The Hilbert function of each generalized Eagon-Northcott module, in the standard graded case, is given in Proposition \tref{P11.4}. Lemma \tref{L10.11} and Corollary \tref{L10.11'} show how to 
compute the Hilbert function of a generalized Eagon-Northcott module in a bi-graded situation.
 The main result of the present section is Proposition \tref{5main}, where  we record the formula for 
$\lambda((S/H)_{(u,s)})$ and $\lambda \((\Cal E_{\pmb a}/\Cal D_{\pmb a})_{(u,s)}\)$ for each eligible tuple $\pmb a$ and each bi-degree $(u,s)$.

\definition{Definition \tnum{D10.5}} Let $P$ be a ring, $E$ and $F$ be free $P$-modules of rank $2$ and $c$, respectively, and $\Psi\:F\to E$ be a homomorphism of  $P$-modules. Define the generalized Eagon-Northcott module $\operatorname{EN}[\Psi,P,r]$ by $$\operatorname{EN}[\Psi,P,r]=\cases
\cok (E^*\t\W^2F\to F)&\text{if $r=-1$}\\
P/I_2(\Psi)&\text{if $r=0$}\\
\Sym_r(\cok \Psi)&\text{if $1\le r$.}\endcases$$
The defining map for $\operatorname{EN}[\Psi,P,-1]$ sends $u\t v$ to $[\Psi^*(u)](v)$. When there is no ambiguity about the ring $P$, we suppress the $P$ and write $\operatorname{EN}[\Psi,r]$ in place of $\operatorname{EN}[\Psi,P,r]$. \enddefinition

\definition{Convention}We define  the binomial
coefficient $\binom{j}{i}$ for all integers   $i$ and $j$  by
$$\binom{j}{i}=\left\{ \matrix\format \c&\quad\l\\
\dfrac{j(j-1)\cdots(j-i+1)}{i!}&\text{if $0<i$,}\\\vspace{5pt} 
1&\text{if $0=i$, and}\\\vspace{5pt} 0&\text{if
$i<0$.}\endmatrix \right.$$
 If $i$ and $j$ are integers with $0\le j$, then $\binom{j}{i}=\binom{j}{j-i}.$
 If $i$ is a nonnegative integer, then $\binom{-1}{i}=(-1)^{i}.$\enddefinition

\proclaim{Proposition \tnum{P11.4}} Let $P$ be a standard graded polynomial ring over a field and let $\psi$ be a $2\times c$ matrix of linear forms in $P$. 
Let $F=P(-1)^{c}$ and $E=P^2$ and view $\psi$ as a map $\psi\:F\to E$.
Assume that $\htt I_2(\psi)=c-1$ and let $D$ be the Krull dimension of $P/I_2(\psi)$. If $r$ and $s$ are integers, with $-1\le r\le c-1$, then
$$\lambda(\operatorname{EN}[\psi,r]_s)=(r+1)\binom{s+D-2}{s}+c\binom{s+D-2}{s-1}.\tag\tnum{ell} $$
\endproclaim
\remark{Remarks}
\item{1.} Notice that both sides of (\tref{ell}) are zero when $s<0$.
\item{2.} If $D=0$, then the right side of (\tref{ell}) is equal to 
$$\cases r+1,&\text{if $s=0$},\\c-(r+1),&\text{if $s=1$,} \\0,&\text{$2\le s$}. \endcases$$
\item{3.} If $D=1$, then the right side of (\tref{ell}) is equal to 
$$\cases r+1,&\text{if $s=0$,}\\c ,&\text{if $1\le s$.} \endcases$$
 \endremark\demo{Proof}The proof is by induction on $D$. 
Start with $D=0$. In this case, the number of variables in $P$ is equal to  $\htt I_2(\psi)=c-1$. In particular, $\lambda(P_1)=c-1$.
 First, fix $r\ge 1$. In this case, $\operatorname{EN}[\psi,r]$ is minimally presented by
$$S_{r-1}E\t F\to S_rE\to \operatorname{EN}[\psi,r]\to 0,$$which is the same as $$P(-1)^{rc}\to P^{r+1}\to \operatorname{EN}[\psi,r]\to 0.$$ It is clear that $\lambda(\operatorname{EN}[\psi,r]_0)=r+1$. One may read that
$$\lambda(\operatorname{EN}[\psi,r]_1)=(r+1)\lambda(P_1)-rc\lambda(P_0)=c-1-r.$$ We know that $I_2(\psi)$ kills $\operatorname{EN}[\psi,r]$. However, $I_2(\psi)$ is equal to the square of the maximal ideal of $P$ (notice that $I_2(\psi)\subseteq \maxm^2$ and both ideals of $P$ are minimally generated by $\binom c 2$ elements of $P_2$), and $\operatorname{EN}[\psi,r]$ is generated in degree zero; so $\operatorname{EN}[\psi,r]_s=0$ for all $s\ge 2$.

It is very easy to see that the assertion is correct for $r=0$.
We now consider $r=-1$. The module  $\operatorname{EN}[\psi,r]$ is minimally presented by
$$E^*\t \tW^2F\to F\to \operatorname{EN}[\psi,r]\to 0,$$which is the same as $$P(-2)^{2\binom{c}2}\to P(-1)^{c}\to \operatorname{EN}[\psi,r]\to 0.$$ We can now read that
$$\lambda(\operatorname{EN}[\psi,r]_s)= \cases 0,&\text{if $s=0$},\\c,&\text{if $s=1$, and}\\c\lambda(P_1)-2\binom{c}2=0,&\text{if $s=2$}.\endcases$$ Once again, all of the generators of $\operatorname{EN}(\psi,r)$ have the same degree. As soon as we know that $\operatorname{EN}[\psi,r]_2=0$, then we know that $\operatorname{EN}[\psi,r]_s=0$ for all $s\ge 2$. 

Now we treat positive $D$. Let $x$ be a linear form in $P$ that is regular on
 $P/I_2(\psi)$. Write $\bar P$ for $P/(x)$ and $\bar \psi$ for $\psi\t _P\bar P$. The module $\operatorname{EN}[r,\psi]$ is perfect (in the sense of \cite{\rref{BH}, Def\.~1.4.14}) and has the same associated primes as $P/I_2(\psi)$. It follows that
$$0\to \operatorname{EN}[\psi,r](-1)@> x>> \operatorname{EN}[\psi,r]\to \operatorname{EN}[\bar\psi,r]\to 0$$ is an exact sequence; and therefore  $\lambda(\operatorname{EN}[\psi,r]_s)=\sum\limits_{i=0}^s\lambda(\operatorname{EN}[\bar \psi,r]_i)$.   
\qed\enddemo

We now study the Hilbert function of the generalized Eagon-Northcott modules in a bi-graded situation. The main algebraic tool is Lemma \tref{L10.11}, which has nothing to do with grading. In Corollary \tref{L10.11'}, we apply Lemma \tref{L10.11} to the bigraded case of interest.

\proclaim{Lemma \tnum{L10.11}} Adopt the notation of Definition \tref{D10.5}. Assume that $F=F'\p F''$ for free modules $F'$ and $F''$ where $F''$ has rank $1$. Let $\Psi'\:F'\to E$ be the restriction of $\Psi$ to $F'$ and $\Psi''\:F''\to E$ be the restriction of $\Psi$ to $F''$. Assume that
$$\grade I_2(\Psi)\ge c-1.$$ If $0\le r\le c-1$, then there is a short exact sequence
$$ 0\to \operatorname{EN}[\Psi',r-1]\t F''@> \iota >> \operatorname{EN}[\Psi',r]@>\pi >> \operatorname{EN}[\Psi,r]\to 0,$$
where $\pi$ is the natural surjection and 
$$\iota(m\t v)=\cases m\cdot \Psi''(v)&\text{for $1\le r$}\\
\W^2\Psi(m\w v)&\text{for $r=0$}.\endcases$$
\endproclaim
\demo{Proof} Recall that the generalized Eagon-Northcott complex that is associated to $\operatorname{EN}[\Psi,r]$ is $\frak E\frak N[\Psi,r]_{\bullet}$ with
$$\frak E\frak N[\Psi,r]_{p}=\cases \Sym_{r-p}E\t \W^pF&\text{if $0\le p\le r$}\\
D_{p-r}E^*\t \W^{p+1}F&\text{if $r+1\le p$}.\endcases$$Recall also, that if $\grade I_2(\Psi)\ge c-1$, then $\frak E\frak N[\Psi,r]_{\bullet}$ is a resolution of $\operatorname{EN}[\Psi,r]$. In the present situation, the decomposition $F=F'\p F''$ induces a short exact sequence of modules
$$0\to \W^pF'\to \W^pF\to \W^{p-1}F'\t F''\to 0,$$ for all $p$. Furthermore, these short exact sequences of modules induce a short exact sequence of complexes
$$0\to \frak E\frak N[\Psi',r]_{\bullet}\to \frak E\frak N[\Psi,r]_{\bullet}\to \frak E\frak N[\Psi',r-1]_{\bullet}[-1]\t F''\to 0,$$ for all $r$. The corresponding long exact sequence of homology includes
$$\Hgy_1(\frak E\frak N[\Psi,r]_{\bullet})\to  \operatorname{EN}[\Psi',r-1]\t F''@> \iota >> \operatorname{EN}[\Psi',r]@>\pi >> \operatorname{EN}[\Psi,r]\to 0.$$ The hypothesis $\grade I_2(\Psi)\ge c-1$ ensures that $\Hgy_1(\frak E\frak N[\Psi,r]_{\bullet})=0$. \qed\enddemo

 \proclaim{Corollary \tnum{L10.11'}} Retain the hypotheses of Lemma \tref{L10.11}. Suppose that the ring $P$ is equal to $P'[x,y]$ where $P'$ is a standard graded polynomial ring over the field $k$ and $x$ and $y$ are new variables.
View $P$ as a bi-graded ring. The variables $x$ and $y$ have degree $(1,0)$. Each variable from $P'$ has degree $(0,1)$. Suppose $\Psi'$ is a $(c-1)\times 2$ matrix of linear forms from $P'$ and $\Psi''=\[\smallmatrix y\\x\endsmallmatrix\]$. Let $R$ be the standard graded polynomial ring $k[x,y]$. If $(u,s)$ is any bi-degree and $r$ is any integer with  
$0\le r\le c$, then
$$\lambda(\operatorname{EN}[\Psi,P,r]_{(u,s)})=\lambda(R_u)\lambda(\operatorname{EN}[\Psi',P',r]_{s})-\lambda(R(-1)_u)\lambda(\operatorname{EN}[\Psi',P',r-1]_{s}).$$
\endproclaim

\demo{Proof}Apply Lemma \tref{L10.11}
to obtain  the short exact sequence 
$$0\to \operatorname{EN}[\Psi',P,r-1](-1,0)\to \operatorname{EN}[\Psi',P,r]\to \operatorname{EN}[\Psi,P,r]\to 0.$$ We have    $P=R\t_k P'$. The map $\Psi'\: P(-1)^{m-2}\to P^{2}$ is the same as $$1\t\Psi'\: R\t_k P'(-1)^{m-2}\to R\t_k{P'}^{2};$$ and therefore,   $\operatorname{EN}[\Psi',P,r]=R\t_k \operatorname{EN}[\Psi',P',r]$. It follows that
$$\align &\lambda((\operatorname{EN}[\Psi,P,r])_{(u,s)})\\&{}= \lambda((R\t_k\operatorname{EN}[\Psi',P',r])_{(u,s)})- \lambda((R(-1)\t_k\operatorname{EN}[\Psi',P',r-1])_{(u,s)})\\&{}= \lambda(R_u)\lambda(\operatorname{EN}[\Psi',P',r]_{s})- \lambda(R(-1)_u)\lambda(\operatorname{EN}[\Psi',P',r-1]_{s}). \qed \endalign$$
\enddemo

 The rest of this section is devoted to proving Proposition \tref{5main}. Adopt the notation of Definition \tref{CN} with (\tref{3.1}). Recall the notion of eligible $k$-tuple $\pmb a=(a_1,\dots,a_k)$, as well as $f(\pmb a)$ and  $r(\pmb a)$, from the statement of Theorem \tref{Kupn}. In \cite{\rref{KPU-d}, Def\.~3.1} we put a total order on the set of eligible tuples. For eligible tuples $\pmb b>\pmb a$ we define ideals $\Cal E_{\pmb b}\subseteq \Cal E_{\pmb a}$ of $A$  by induction. There is no convenient way to denote the eligible tuple which is immediately larger than a particular eligible tuple $\pmb a$; consequently, we define two parallel collections of ideals $\{\Cal E_{\pmb a}\}$ and $\{\Cal D_{\pmb a}\}$ simultaneously. The ideal $\Cal D_{\emptyset}$ is equal to zero. If $\pmb a$ is an eligible tuple of positive length, then $\Cal D_{\pmb a}=\sum\limits_{\pmb b>\pmb a} \Cal E_{\pmb b}$.  If $\pmb a$ is an arbitrary eligible tuple, then
$$\Cal E_{\pmb a}=\Cal D_{\pmb a}+T^{\pmb a}T_{k+1,1}^{f(\pmb a)}(T_{k+1,1},\dots,T_{k+1,r(\pmb a)}).$$
We have a filtration of $K^{(n)}$: 
$$(0)\subsetneq \Cal E_{\emptyset}\subsetneq \dots \subsetneq \Cal E_{0^{\ell-1}}=K^{(n)},$$where $0^s$ is the $s$-tuple $(0,\dots,0)$. It is also shown in \cite{\rref{KPU-d}, 
Thm\.~3.17} that the factor module $\Cal E_{\pmb a}/\Cal D_{\pmb a}$ is isomorphic to the generalized Eagon-Northcott module 
$$\operatorname{EN}[\psi_{>k},S/P_{k},r(\pmb a)-1](-t_{\pmb a}),$$
where $\psi_{>k}$ is the submatrix $\bmatrix \psi_{k+1}&\cdots&\psi_{\ell}\endbmatrix$ of $\psi$, $P_{k}$ is the ideal $$P_k=\(\{T_{i,j}\mid 1\le i\le k, 1\le j\le \sigma_i+1\}\)$$ of $S$, and $t_{\pmb a}$ is the twist 
$$t_{\pmb a}=\cases (0,\sum\limits_{u=1}^k a_u+f(\pmb a)+1),&\text{if $k<\rho$,}\\
(f(\pmb a)+1,\sum\limits_{u=1}^k a_u),&\text{if $k=\rho$}.\endcases$$

\proclaim{Proposition \tnum{5main}}
    Adopt the notation of Definition \tref{CN} with $(\tref{3.1})$. Let $(u,s)$ be an arbitrary bi-degree. 
\medskip\flushpar{\bf(a)}  
$$\lambda((S/H)_{(u,s)})=\lambda(R_u)
\(\binom{s+1}{s}+(m-2)\binom{s+1}{s-1}\)
-\lambda(R(-1)_u)(m-2)\binom{s+1}{s-1}$$
 
\medskip\flushpar{\bf(b)} 
$$\lambda\(\(\Cal E_{\emptyset}/\Cal D_{\emptyset}\)_{(u,s)}\)=\left\{\matrix
\lambda(R_u)\[r(\emptyset)\binom{s-f(\emptyset)}{s-f(\emptyset)-1}+(m-2)\binom{s-f(\emptyset)}{s-f(\emptyset)-2}\]\hfill\\
-\lambda(R(-1)_u)\[(r(\emptyset)-1)\binom{s-f(\emptyset)}{s-f(\emptyset)-1}+(m-2)\binom{s-f(\emptyset)}{s-f(\emptyset)-2}\]\hfill\endmatrix\right.
$$
 \medskip\flushpar{\bf(c)} If $\pmb a=(a_1)$ is  an eligible $1$-tuple, then $\lambda\(\(\Cal E_{\pmb a}/\Cal D_{\pmb a}\)_{(u,s)}\)$ is equal to 
$$\cases \chi(a_1\le s)\lambda(R(a_1\sigma_1-n)_u),&\text{if $\rho=1$, or}\\
\chi(a_1+f(a_1)+1\le s)\left(\matrix \phantom{+}\lambda(R_u)(a_1\sigma_1-n+1+\sigma_2(s-a_1))\hfill\\-\lambda(R(-1)_u)(a_1\sigma_1-n+\sigma_2(s-a_1))\hfill\endmatrix\right),
&\text{if $\rho=2$.}\endcases$$
\medskip\flushpar{\bf(d)} If  $\pmb a=(a_1,a_2)$ is  an eligible $2$-tuple then
$$\lambda\(\(\Cal E_{\pmb a}/\Cal D_{\pmb a}\)_{(u,s)}\)=\chi(s=a_1+a_2)\lambda(R(
a_1\sigma_1+a_2\sigma_2-n)_u).$$
\endproclaim
  \demo{Proof}
For (a) and (b) we apply Corollary \tref{L10.11'} with $P'=k[T_1,\dots,T_m]$ and $\Psi'$ equal to the first $m-2$ columns of $\psi$. Thus,
$$\align \lambda((S/H)_{(u,s)})&{}= \lambda((\operatorname{EN}[\psi,S,0])_{(u,s)})
 \\
&{}=\lambda(R_u)\lambda(\operatorname{EN}[\Psi',P',0]_{s})- \lambda(R(-1)_u)\lambda(\operatorname{EN}[\Psi',P',-1]_{s})\endalign$$
and 
$$\align \lambda\((\Cal E_{\emptyset}/\Cal D_{\emptyset})_{(u,s)}\)&{}=\lambda\((\operatorname{EN}[\psi,S,r(\emptyset)-1](0,-f(\emptyset)-1))_{(u,s)}\)  \\
&{}=\left\{\matrix \phantom{-}\lambda(R_{u})
\lambda\((\operatorname{EN}[\Psi',P',r(\emptyset)-1](-f(\emptyset)-1))_s\)\hfill\\
-\lambda(R(-1)_{u})\lambda\((\operatorname{EN}[\Psi',P',r(\emptyset)-2](-f(\emptyset)-1))_s\).\hfill\endmatrix\right.\endalign$$
Apply Proposition \tref{P11.4}, with $c=m-2$ and $D=3$,  to establish (a) and (b). 

Take $\pmb a=(a_1)$ to be an eligible $1$-tuple with $\rho=1$. Apply Corollary \tref{L10.11'} with $P'=k[T_m]$ and $\Psi'$ equal to the zero map. In this case, $r(\pmb a)=1$, $f(\pmb a)+1=n-a_1\sigma_1$, $\operatorname{EN}[0,P',0]=P'$, and $\operatorname{EN}[0,P',-1]=0$. We have 
$$\align \lambda\((\Cal E_{\pmb a}/\Cal D_{\pmb a})_{(u,s)}\)&{}=\lambda\(\(\operatorname{EN}\[\[\smallmatrix y\\x\endsmallmatrix\],P'[x,y],0\](a_1\sigma_1-n,-a_1)\)_{(u,s)}\)  \\
&{}=\lambda(R(a_1\sigma_1-n)_u)\lambda(P'(-a_1)_s).\endalign$$
If $\pmb a=(a_1,a_2)$ is an eligible $2$-tuple, then $\rho$ must equal $2$, $r(\pmb a)=1$, $f(\pmb a)+1=n-a_1\sigma_1-a_2\sigma_2$, 
$$\align \lambda\((\Cal E_{\pmb a}/\Cal D_{\pmb a})_{(u,s)}\)&{}=\lambda\(\(\operatorname{EN}\[\[\smallmatrix y\\x\endsmallmatrix\],k[x,y],0\](a_1\sigma_1+a_2\sigma_2-n,-a_1-a_2)\)_{(u,s)}\)\\&{}=\lambda(R(a_1\sigma_1+a_2\sigma_2-n)_u)\lambda(k(-a_1-a_2)_s).\endalign$$Finally, let $\pmb a=(a_1)$ be an eligible $1$-tuple with $\rho=2$. Apply Corollary \tref{L10.11'} with $P'=k[T_{2,1},\dots,T_{2,\sigma_2+1}]$ and $\Psi'=\psi_2$ to see that
$$\align \lambda\((\Cal E_{\pmb a}/\Cal D_{\pmb a})_{(u,s)}\)&{}=\lambda\(\(\operatorname{EN}\[\[\matrix \psi_2&\psi_3\endmatrix\],P'[x,y],r(\pmb a)-1\](0,-a_1-f(\pmb a)-1)\)_{(u,s)}\)\\&{}=\left\{\matrix \phantom{-}\lambda(R_{u})
\lambda\((\operatorname{EN}[\psi_2,P',r(\pmb a)-1](-a_1-f(\pmb a)-1))_s\)\hfill\\
-\lambda(R(-1)_{u})\lambda\((\operatorname{EN}[\psi_2,P',r(\pmb a)-2](-a_1-f(\pmb a)-1))_s\).\hfill\endmatrix\right.\endalign$$
Apply Proposition \tref{P11.4}, with $c=\sigma_2$ and $D=2$,  to complete the calculation.
\qed
  \enddemo

 \bigskip
 \bigskip
\bigpagebreak
\SectionNumber=6\tNumber=1
\flushpar{\bf \number\SectionNumber.\quad  The resolution of $I^s$.}

We resolve every power of the ideal $I$ of Definition \tref{CN}. Our answer is expressed in terms of the parameter ``$a$'', which is equal to the number of non-linear columns in the matrix which presents $I^s$.  The resolution  depends on the shape of the partition $\pmb \sigma$ which corresponds to $I$. 

\proclaim{Theorem \tnum{pwrs4}} Let $I$ be the ideal of Definition \tref{CN} and $s$ be a positive integer. The minimal homogeneous resolution of $I^s$ has the form
$$0\to  R(-sd-1)^{b} \p \F 
\to R(-sd)^{b_0}\to I^s\to 0,$$
with $b_0=b+a+1$. 
\medskip\flushpar{\bf(1)} If $\rho=1$, then
$\F= \sum\limits_{u=0}^{a-1} R(-sd+u\sigma_1-n)$, $b=sd+\binom a2\sigma_1-an$, and
 $$\tsize a=
\min\left\{s,\left\lceil\frac{n-1}{\sigma_1} \right \rceil\right\}.$$  

\medskip\flushpar{\bf(2)}
If $\rho=2$ and $\sigma_1>\sigma_2$, then $\F=
\sum\limits_{u=0}^{a-1} R(-sd+u(\sigma_1-\sigma_2)+(s-1)\sigma_2-n)$,
 $$\matrix   
b=s(d+a\sigma_2)+\binom a2(\sigma_1-\sigma_2)-a(n+\sigma_2),\text{ and}\hfill \\ \vspace{3pt}
a=\cases \min\left\{s,\left\lceil\frac{n-(s-1)\sigma_2-1}{\sigma_1-\sigma_2} \right \rceil\right\},&\text{if $s\le \frac{n-2}{\sigma_2}+1$,}\\
0,&\text{if $\frac{n-1}{\sigma_2}+1\le s$.}\\
\endcases\hfill \endmatrix $$

\medskip\flushpar{\bf(3)}
If $\rho=2$ and $\sigma_1=\sigma_2$,  then
$\F=
 R(-sd +(s-1)\sigma_2-n)^a,$$$b=s(d+a\sigma_2)-a(n+\sigma_2)\quad\text{and}\quad 
a=\cases 
s,&\text{if $s\le \frac{n-2}{\sigma_2}+1$ and}\\
0,&\text{if $\frac{n-1}{\sigma_2}+1\le s$.}
\endcases $$ 
\endproclaim
\remark{Remark} It is worth noting that the non-linear columns in the presenting matrix for $I^s$ all have the same degree for $\sigma_1=\sigma_2$;  however, these columns
have distinct degrees in the other two cases.\endremark 

\demo{Proof}
The  ring $S$ is   bi-graded  and the   quotient map $$S\onto S/\Cal A=\Cal R(I)$$ sends $S_{(u,s)}\onto R_uI^st^s=I^s_{u+sd}t^s$, where $d$ is the degree of the generators of $I$; so, 
$$\lambda(I^s_{u+ds})=\lambda((S/\Cal A)_{(u,s)}),$$ and, for all integers $s$ and $\frak z$, 
$$\lambda(I^s_{\frak z})=\lambda((S/\Cal A)_{(\frak z-ds,s)}).$$ The short exact sequence $$0\to\Cal A/H\to  S/H\to S/\Cal A\to 0$$ gives
$$\lambda((S/\Cal A)_z)=\lambda((S/H)_z)-\lambda((\Cal A/H)_z).$$
Write $\overline{\phantom{x}}$ to mean image in $A$, as in the proof of Theorem \tref{main1}. The element $\bar g/\bar y^n$ of the quotient field of $A=S/H$ has degree $(0,1)$, since $\bar g$ has degree $(n,1)$ and $\bar y^n$ has degree $(n,0)$, and the isomorphism $\bar g/\bar y^n\: K^{(n)}\to \Cal AA$ of ideals satisfies $\lambda(K^{(n)}_{(u,s-1)})=\lambda((\Cal A/H)_{(u,s)})$. It follows that
$$\lambda(I^s_{\frak z})=\lambda((S/H)_{(\frak z-ds,s)})-\lambda(K^{(n)}_{(\frak z-ds,s-1)}).$$  
We have identified a filtration
$$\{\Cal E_{\pmb a}\mid \text{ $\pmb a$ is an eligible tuple}\}$$ of $K^{(n)}$; thus
$$\lambda(I^s_{\frak z})=\lambda((S/H)_{(\frak z-ds,s)})-\sum_{\pmb a}\lambda\((\Cal E_{\pmb a}/\Cal D_{\pmb a})_{(\frak z-ds,s-1)}\).\tag\tnum{lam4}$$ 
 Each length on the right hand side of  (\tref{lam4}) has been calculated in Proposition \tref{5main}. We have  
$$\lambda(I^s_{\frak z})=   b_0\lambda(R(-sd)_{\frak z})-N_1\lambda(R(-sd-1)_{\frak z})-N_2\tag\tnum{lam14}$$   for 
$$N_2=\cases \sum\limits_{(a_1,a_2)\text{ eligible}} \chi(a_1+a_2=s-1)\lambda(R(-sd+a_1\sigma_1+a_2\sigma_2-n)_{\frak z}),&\text{if $\rho=2$,}\\
\sum\limits_{(a_1)\text{ eligible}}\chi(a_1\le s-1)\lambda(R(-sd+a_1\sigma_1-n)_{\frak z}),&\text{if $\rho=1$},\\
\endcases$$
and integers $b_0$ and $N_1$. (There is no difficulty in recording the exact values of $b_0$ and $N_1$, but this is not necessary.)

When $\rho=2$, we simplify $N_2$ by replacing $a_2$ with $s-1-a_1$. The parameter $a_1$ must satisfy:
$$0\le a_1\le s-1\and a_1\sigma_1+(s-1-a_1)\sigma_2< n.$$Thus, 
$$N_2=N_2'\lambda(R(-sd-1)_{\frak z})+N_2''\tag\tnum{form'}$$ for 
$
N_2'=\sum\limits_{a_1=0}^{s-1}\chi(a_1(\sigma_1-\sigma_2)+(s-1)\sigma_2-n= -1)$ and
$$\eightpoint N_2''=\sum\limits_{a_1=0}^{s-1} \chi(a_1(\sigma_1-\sigma_2)+(s-1)\sigma_2-n\le -2)\lambda(R(-sd+a_1(\sigma_1-\sigma_2)+(s-1)\sigma_2-n)_{\frak z}).$$ 
 When $\rho=1$, we write $N_2$ in the form (\tref{form'}) with
$
N_2'=\sum\limits_{a_1=0}^{s-1}\chi(a_1\sigma_1-n= -1)$ and
$$
N_2''=\sum\limits_{a_1=0}^{s-1} \chi(a_1\sigma_1-n\le -2)\lambda(R(-sd+a_1\sigma_1-n)_{\frak z}). $$
  Let $b=N_1+N_2'$.  Apply Lemma \tref{L19.14} to see that the minimal resolution of $I^s$ is
$$0\to  R(-sd-1)^{b} \p  
\F 
 \to R(-sd)^{b_0}\to I^s\to 0,$$
for $\F$ equal to  $$\sum\limits_{a_1=0}^{s-1} \chi(a_1(\sigma_1-\sigma_2)+(s-1)\sigma_2-n\le -2)R(-sd+a_1(\sigma_1-\sigma_2)+(s-1)\sigma_2-n),$$ if $\rho=2$; or
$$\sum\limits_{a_1=0}^{s-1} \chi(a_1\sigma_1-n\le -2)R(-sd+a_1\sigma_1-n),$$
 if $\rho=1$ . 
Notice that the rank of $\F$ is equal to the number of non-linear columns in the presenting matrix for $I^s$.
We next express $\F$ in a more transparent manner.

When $\rho=1$, the constraint $ a_1\sigma_1-n\le -2$ is equivalent to 
$$a_1\le \left\lfloor \frac{n-2}{\sigma_1}\right\rfloor=\left\lceil \frac{n-1}{\sigma_1}\right\rceil-1$$
and $$\F=\sum\limits_{a_1=0}^{a-1} R(-sd+a_1 \sigma_1 -n),$$ for $a=\min\{s,\left\lceil \frac{n-1}{\sigma_1}\right\rceil\}$.  

 Take $\rho=2$. The parameter $a_1$ is non-negative; so, $\F$ is zero if $\frac{n-1}{\sigma_2}+1\le s$. We think about $s\le \frac{n-2}{\sigma_2}+1$. If $\sigma_2=\sigma_1$, then $\chi((s-1)\sigma_2-n\le -2)=1$ and $$\F=\sum\limits_{a_1=0}^{s-1} R(-sd+(s-1)\sigma_2-n)=R(-sd+(s-1)\sigma_2-n)^s.$$ If $\sigma_1>\sigma_2$, then 
$$\align a_1(\sigma_1-\sigma_2)&{}+(s-1)\sigma_2-n\le -2\iff {}\\ &a_1\le \left\lfloor \frac{n-(s-1)\sigma_2-2}{\sigma_1-\sigma_2}\right\rfloor= \left\lceil\frac{n-(s-1)\sigma_2-1}{\sigma_1-\sigma_2} \right \rceil-1,\endalign$$
and 
$\F=\sum\limits_{a_1=0}^{a-1} R(-sd+a_1(\sigma_1-\sigma_2)+(s-1)\sigma_2-n)$, for
$$a=\cases \min\left\{s,\left\lceil\frac{n-(s-1)\sigma_2-1}{\sigma_1-\sigma_2} \right \rceil\right\},&\text{if $s\le \frac{n-2}{\sigma_2}+1$,}\\
0,&\text{if    $\frac{n-1}{\sigma_2}+1\le s$.}\\
\endcases $$  

Finally, we see  that the values of $b_0$ and $b$ are completely determined by $a$. Indeed, rank is additive on short exact sequences; so,
$b_0=b+a+1$. Also, $I^s$ is generated by the maximal minors of the matrix which presents $I^s$. In other words, $sd$ is equal to the sum of the column degrees of this presenting matrix; that is,
$$\matrix b=sd+\sum\limits_{a_1=0}^{a-1}(a_1(\sigma_1-\sigma_2)+(s-1)\sigma_2-n)\hfill \\\phantom{b}=sd+\binom a2(\sigma_1-\sigma_2)+(s-1)\sigma_2 a-na. \qed\hfill \endmatrix\tag\tnum{tag4}$$     
\enddemo

\proclaim{Lemma \tnum{L19.14}}Let $M$ be a homogeneous module of projective dimension one over the standard graded polynomial ring $R$. Suppose that all of the generators of $M$ have degree $D$. Suppose further that $b_0$, $b_1$ and  $t_1\le t_2\le \dots\le t_{b_1}$ are integers which satisfy $D<t_1$ and 
$$\lambda(M_{\frak z})=b_0\lambda(R(-D)_{\frak z})-\sum\limits_{i=1}^{b_1}\lambda(R(-t_i)_{\frak z})$$ for all integers $\frak z$. Then the minimal homogeneous resolution of $M$ has the form
$$0\to \bigoplus\limits_{i=1}^{b_1}R(-t_i)\to R(-D)^{b_0}\to M\to 0.$$\endproclaim
\demo{Proof}The hypotheses ensure that the minimal homogeneous resolution of $M$ has the form  
$$0\to \bigoplus\limits_{i=1}^{b_1'}R(-t_i')\to R(-D)^{b_0'}\to M\to 0 \tag\tnum{RR194} $$
for some integers $b_0'$, $b_1'$, and $t_1'\le t_2'\le \dots\le t_{b_1'}'$ with $D<t_1'$.
Use (\tref{RR194}) to compute the Hilbert function of $M$; so
$$b_0\lambda(R(-D)_{\frak z})-\sum\limits_{i=1}^{b_1}\lambda(R(-t_i)_{\frak z})=b_0'\lambda(R(-D)_{\frak z})-\sum\limits_{i=1}^{b_1'}\lambda(R(-t_i')_{\frak z}),$$  for all integers $\frak z$.
 It follows  that the free modules
$$\F=R(-D)^{b_0}\p\bigoplus_{i=1}^{b_1'}R(-t_i')\and \F'=R(-D)^{b_0'}\p\bigoplus_{i=1}^{b_1}R(-t_i)$$ have the same Hilbert function. This forces the free $R$-modules  $\F$ and $\F'$ to be equal; in other words, they have the exact same twists: $b_0=b_0'$, $b_1=b_1'$, and $t_i=t_i'$ for all $i$. \qed
\enddemo

The first two assertions of the following result may be read from the resolution of Theorem \tref{pwrs4}. A different proof of these results may be found  in Corollary \tref{C2.7}. 

\proclaim{Corollary \tnum{reg}} Let $I$ be the ideal of Definition \tref{CN} and $s$ be a positive integer.
\flushpar{\bf (1)} If $\rho=1$, then   $\reg I^s= sd+n-1$ for all $s\ge 1$.
\medskip \flushpar{\bf (2)} If $\rho=2$, then  $\reg I^s= sd$ if and only if  $\frac{n-1}{\sigma_2}+1\le s$.
\medskip\flushpar{\bf(3)} The following statements are equivalent:
\itemitem{\rm (a)} $I^s=(x,y)^{sd}$,
\itemitem{\rm (b)} the minimal homogeneous resolution of $I^s$ has the form 
$$0\to R(-sd-1)^{b-1}\to R(-sd)^{b}\to I^s\to 0,$$ for some $b$, 
\itemitem{\rm(c)} $\rho=2$   and $\frac{n-1}{\sigma_2}+1\le s$, or $\rho=1$ and $n=1$, and 
\itemitem{\rm(d)} $a=0$.
\itemitem{\rm(e)} $\reg I^s=sd$.
\endproclaim
\demo{Proof} We prove (3). The trick (\tref{tag4}) shows that (a) and (b) are equivalent. 
The parameter $a$ is equal to the number of non-linear columns in the presenting matrix for $I^s$, so (d) and (b) are equivalent. The equivalence of (d) and (c) may be read from Theorem \tref{pwrs4}. Assertions (1) and (2) show that (c) and (e) are equivalent.\qed
\enddemo

Let $B$ be a standard graded algebra over a field and let $q_B(s)$ be the Hilbert polynomial of $B$. It follows that $q_B(s)=\lambda(B_s)$ for all large $s$. The postulation number of $B$ is
$$p(B)=\max\{s\mid q_B(s)\neq \lambda(B_s)\}.$$
\proclaim{Corollary \tnum{last}}If $I$ is given in Definition \tref{CN}, then 
$$p(\Cal F(I))=\cases
\lceil \frac{n-1}{\sigma_{2}}\rceil,&\text{if $\rho=2$,}\\
\lceil \frac{n-1}{\sigma_{1}}\rceil -1,&\text{if $\rho=1$.}\endcases$$
\endproclaim
\demo{Proof} The Hilbert function and the Hilbert polynomial of $\Cal F(I)$ may be read from Theorem \tref{pwrs4}: $\Hgy_{\Cal F(I)}(s)$ is equal to ``$b_0$'', written  as a function of $s$ and
$$q_{\Cal F(I)}(s)=\cases sd +\binom a2\sigma_1-an+a+1,&\text{if $\rho=1$,}\\sd+1,&\text{if $\rho=2$,}\endcases$$for $a=\lceil \frac{n-1}{\sigma_1}\rceil$. The calculation of $p(\Cal F(I))$ when $\rho=2$ is explicitly given in Corollary \tref{reg}. A similar calculation is used when $\rho=1$. \qed\enddemo

\medskip

\flushpar {\bf Acknowledgment.}
This work was conducted while the first author was on sabbatical at Purdue University and later was a Visiting Professor at the University of Notre Dame. He appreciates the sabbatical from the University of South Carolina and the hospitality he received at Purdue University and the University of Notre Dame.  Also, we appreciate that Wolmer Vasconcelos made us aware of Conjecture 4.5 in \cite{\rref{HSV}}. This conjecture, later established by \cite{\rref{CHW}}, is the starting point of this project.

\enddocument